\documentclass[preprint,authoryear,3p,number]{elsarticle}

\usepackage{lipsum}
\usepackage{gensymb}
\usepackage{graphicx,psfrag}
\usepackage{enumerate}
\usepackage{mathrsfs}
\usepackage{amsmath,amsfonts,amssymb}
\usepackage{enumerate}
\usepackage{hyperref}
\usepackage{float}
\usepackage{ulem}
\usepackage{color}
\usepackage{natbib}
\setcitestyle{numbers}
\biboptions{sort&compress}

\newcommand{\p}{\partial}
\newcommand{\femn}{$\textrm{FEM}_N$ }
\newcommand{\sn}{$S_N$ }
\newcommand{\fpn}{$\textit{\textrm{FP}}_N$ }

\title{A Finite Element Method for Angular Discretization of the Radiation Transport Equation on Spherical Geodesic Grids}

\author[1,2]{Maitraya K Bhattacharyya\corref{cor1}}
\ead{mbb6217@psu.edu}
\author[1,2,3]{David Radice\corref{cor2}}
\ead{dur566@psu.edu}

\cortext[cor1]{Corresponding author}
\cortext[cor2]{Alfred P.~Sloan Fellow}

\address[1]{Institute for Gravitation and the Cosmos, The Pennsylvania State University, University Park, PA 16802, USA}
\address[2]{Department of Physics, The Pennsylvania State University, University Park, PA 16802, USA}
\address[3]{Department of Astronomy \& Astrophysics, The Pennsylvania State University, University Park, PA 16802, USA}
\begin{document}
\begin{abstract}
	Discrete ordinate~($\textit{\textrm{S}}_{N}$) and filtered spherical harmonics~($\textit{\textrm{FP}}_{N}$) based schemes have been proven to be robust and accurate in solving the Boltzmann transport equation but they have their own strengths and weaknesses in different physical scenarios. We present a new method based on a finite element approach in angle that combines the strengths of both methods and mitigates their disadvantages. The angular variables are specified on a spherical geodesic grid with functions on the sphere being represented using a finite element basis. A positivity-preserving limiting strategy is employed to prevent non-physical values from appearing in the solutions. The resulting method is then compared with both~\sn and~\fpn schemes using four test problems and is found to perform well when one of the other methods fail.
\end{abstract}


\maketitle

\section{Introduction}
The transport of radiation carriers like leptons, photon, hadrons and ions in different types of media is governed by the Boltzmann equation, which has applications in many areas of science and engineering. For instance, in the field of biophotonics and biomedical imaging, the Boltzmann equation is used to describe photon transport for X-ray and tomography experiments~\cite{WanTuc2013,ChaBouKim2008}. In radiology, the transport of photons is used to estimate dosages of nuclear medicine in a clinical setting~\cite{Bed2019}. Similarly, in the field of metrology and electronic circuit design, this equation describes phonons and is used to model the transport of heat~\cite{WheShaTam2014}.  It also proves to be a vital tool in the field of electron microscopy where it describes the scattering of electrons~\cite{FatRez1982}. The Boltzmann equation has many other diverse applications from weather and climate modeling~\cite{ThoSta2002, CloSheMla2005} to the design and characterization of nuclear detectors~\cite{WagPepDou2011}. In the field of astrophysics, these equations are especially important in multimessanger astrophysics to describe neutrinos which play an important role in driving the winds in remnants of binary neutron star mergers~\cite{DesOttBur2008, PerRosCab2014,FujSekKiu2017}, determining the composition of ejecta in such mergers~\cite{RadGalLip2016,PerRadBer2017,FouHaaDue2016,SekKiuKyu2016}, core-collapse supernova~\cite{MezEndMes2020} and gamma-ray bursts~\cite{PerYasArc2017} among others. With increasing instrument sensitivities, accurate modeling of astrophysical observables demand better models of neutrino transport. In binary neutron star simulations, for instance, several efforts have been made to move away from phenomenological models~\cite{SekKiuKyu2011,BetWil1985} to more realistic ones. One of the main roadblocks towards more accurate radiation transport models is computational cost: the Boltzmann transport equation, which describes the time evolution of the distribution function of radiation carriers is a function of seven independent variables. In many realistic scenarios this equation has to be solved without any symmetry considerations for each particle species. This poses a significant computational challenge which cannot be tackled by a simple-minded brute-force approach.

Several approaches have been adopted to solve the transport equation. These can be broadly classified into two categories: approximate methods and methods which solve the full Boltzmann equation. In the case of approximate methods, the Boltzmann equation is substituted by a more manageable approximate equation. Moment based methods, which falls into this category, have been used quite successfully in core-collapse supernova~\cite{Con2015,KurTakKot2016,ConCou2018,RobOttHaa2016,SkiDolBur2019,GlaJusJan2019,RahJusJan2019,LaiEndChu2021} and binary neutron star merger~\cite{FouOCoRob2015,FouHaaDue2016,FouOCoRob2016,RadBerPer2022} simulations. In these schemes, the transport equations are rewritten as a sum of moments of the radiation distribution function up to a particular order. The resulting system is then closed by specifying a closure condition. The one-moment method evolves the zeroth moment of the distribution function, that is the radiation energy density and uses a closure relation to compute the first moment or the radiation flux. Similarly, in the two-moment method, the zeroth and the first moments are evolved and a closure is used to evaluate the second moment: the radiation pressure. The advantage of moment-based methods is the significantly lower computational cost compared to solving the full transport equation. This advantage, however comes at the price of accuracy, which depends on the order of truncation of the scheme and the choice of closure~\cite{Ric2020}.

The full Boltzmann equation has been solved using both probabilistic and deterministic techniques. The Monte Carlo method~\cite{FleCum1971,FleCan1984,DenUrbEva2007,AbdErnBur2012,RicKasOCon2015}, which uses pseudo-random number generators to directly simulate the transport of radiation carriers is one of the most accurate methods to solve the Boltzmann equation. This method has the advantage of being easily adaptable in both simple and complex geometries and of handling anisotropies in the problem with relative ease. Monte Carlo algorithms are highly parallelizable and have been shown to exhibit excellent scaling capabilities~\cite{AbdErnBur2012}. The downside of this method is the high computational cost required to obtain solutions of sufficiently high accuracy since the stochastic nature of the method introduces statistical noise in solutions. According to the central limit theorem, this noise scales as~$N^{-1/2}$, where~$N$ is the number of particles~\cite{AbdErnBur2012}. Furthermore, explicit Monte Carlo schemes are computationally expensive at high optical depth~\cite{CleGen2014,PoeValBer2020} owing to the fact that a large number of particle interactions have to be considered in scattering dominated regions. Implicit schemes~\cite{PoeVal2020,SteHei2022,SteHei2022b}, which employ acceleration techniques, like discrete diffusion Monte Carlo schemes~\cite{CleGen2014} have been proposed to improve computational efficiency.  However, these methods are generally problematic to implement in the context of relativistic radiation hydrodynamics since the relativistic diffusion equation is not well-posed.

Out of the deterministic methods, we focus our attention on the discrete ordinate method or~\sn and the filtered spherical harmonics method or~$\textrm{\textit{FP}}_N$. Both have their merits and demerits, but are known to provide robust, accurate solutions for particular problems. In the ~\sn method~\cite{MihWei1984,GodLiu2012,LarMor2010,Liv2003,ChaMul2020}, the distribution function for radiation carriers is evolved directly by discretizing it along~$N$ angular bins. Modern numerical implementations of the resulting equations have good computational efficiency and can produce solutions in the optically thick limit with high accuracy. This method is also very efficient in handling highly-directional beams of radiation. The main problem with this approach is the absence of rotational invariance, which gives rise to ``ray effects". These are artifacts of the angular discretization which appear as oscillations in the spatial domain~\cite{Ten2014} and heavily diminish the quality of the solution, especially in regions of low scattering. These effects can be reduced by increasing the number of discrete ordinates and employing efficient filtering strategies~\cite{HauHen2019}. 

The spherical harmonics method or~$\textrm{\textit{P}}_N$~\cite{McCHau2010,RadAbdRez2013,McCEvaLow2008,McCHolBru2008} on the other hand, expands the distribution function in a basis of spherical harmonics. This choice of basis ensures explicit rotational invariance. The main disadvantage of this method is the appearance of non-physical oscillations in regions of low optical depth, which can cause the distribution function to acquire negative values. When compared to the~\sn method, the~$\textrm{\textit{P}}_N$ method performs poorly with rays or beams of radiation. A modified approach, called~\fpn uses filtering techniques to remove the Gibbs' oscillations from the system~\cite{McCHau2010,RadAbdRez2013}. The strength of filters is usually controlled by a tunable parameter called the effective opacity. A filter with large effective opacity may be successful in eliminating negative values from the distribution function, but at the cost of a reduced solution quality. On the other hand, choosing a very low value for the effective opacity may not eliminate all oscillations or non-physical values from the solution. The choice of this parameter is also problem specific. 

In this paper, we introduce a ``best of both worlds" approach using a finite element method in angle:~$\textrm{FEM}_{N}$. Finite element approaches for angular discretization of the Boltzmann equation have been proposed to reduce ``ray artifacts'' seen in the~\sn method~\cite{MorWarLow2003,DulBarPri2014}. Several discretization strategies have been proposed~\cite{CopLapRav1990,Kan2009,EggSch2016,WanAbeMud2020,Jar2010,GanSin2022}, including the use of wavelets for angular mesh refinement~\cite{BucPaiEat2005}, discontinuous finite element approaches in angle~\cite{KopLat2015} and multi-$P_N$ schemes~\cite{GhaAbbZol2019,FalSciMas2022} which introduce the benefits of the~$P_N$ method within a solid angle. However, most of these strategies discretize the angular variables on a latitude-longitude grid, so special treatment of the poles is necessary where co-ordinate singularities arise. Moreover, positivity preservation in~$P_N$ methods are not straightforward, as will be demonstrated later. Discontinuous finite element methods~\cite{KopLat2015} reduce to the~\sn method at the low resolutions and facilitates the use of efficient sweeping algorithms with local angular refinement. However, methods using a discontinuous finite element basis requires more angular resolution to resolve smooth solutions than those using continuous basis functions. In this paper, we propose a new finite element method in angle where the angular discretization is performed using spherical geodesic grids~\cite{Gir1997} and positivity preserving is ensured using limiters adapted from~\fpn schemes~\cite{LaiHau2019}. The distribution function is expanded in angle in terms of finite element~(FEM) basis functions on the geodesic grid. Every point on the geodesic grid is expressed in Cartesian coordinates, thereby side-stepping the pathologies at the poles. Furthermore, the triangular elements of the grid approximately occupy the same area and mesh refinement strategies in angle are straightforward. In our formulation, the Boltzmann equation can be written down as a system of advection equation with a specific choice of basis functions giving rise to the~\femn method or the~\sn method. Moreover, positivity of the distribution function and the point-wise conservation of the radiation energy density at every time step is ensured by the use of a parameter-free limiter strategy, which is a potential advantage over the filtering strategies used in the~$\textrm{P}_N$ method.

The Boltzmann equation has distinctly different behavior in the high and low optical limits. In regions of high opacity, the equation becomes diffusive while at low opacities, the equations are hyperbolic. A numerical method able to correctly approximate the Boltzmann equation at all optical depths is preferred. This is achieved by using an asymptotic preserving~(AP) discontinuous Galerkin scheme~\cite{LowMor2002,McCEvaLow2008,RadAbdRez2013}. A second order Runge-Kutta method is used for the time stepping when the source terms are not stiff, as is the case with all numerical experiments in this paper. A semi-implicit method~\cite{McCEvaLow2008} may be used for problems with stiff sources.

The paper is organized as follows. In section~\ref{sec:boltzmann}, we introduce basic terminology and the form of the Boltzmann equation with sources to be used in the rest of the paper. For the sake of simplicity, we consider the distribution function to be independent of the frequency of radiation. The description of the~\femn numerical scheme is provided in section~\ref{sec:scheme}. The geodesic grid used for angular discretization and strategies for generation and refinement are described in sub-section~\ref{subsec:angular}. The specific ansatz taken for the distribution function is used to arrive at a general set of coupled advection equations with specific choices of basis recovering the~$\textrm{\textit{S}}_N$,~\fpn and~\femn schemes. A description of the spatial discretization using an asymptotic preserving discontinuous Galerkin scheme is provided in sub-section~\ref{subsection:dg}. The time integrator and the~\femn limiter for positivity preservation and the filter for~\fpn methods are discussed in sub-sections~\ref{subsec:time} and~\ref{subsec:limfilt} respectively. Finally, in section~\ref{section:results}, we perform a systematic comparison between the three methods and discuss the advantages and disadvantages of the~\femn method over its counterparts. Section~\ref{sec:conclusions} summarizes the main results of the paper.

\section{The Boltzmann equation} \label{sec:boltzmann}
Radiation carriers, which are point particles carrying energy~$\epsilon = h \nu$, can be described by a distribution function~$F$, which in its full generality is a function of seven variables: time, three spatial coordinates~$x^i$, two angular coordinates~$\Omega = (\phi, \theta)$ representing the direction of propagation, and the frequency of radiation~$\nu$. The distribution function is defined such that
\begin{align}
	dN = \frac{h^3 \nu^2}{c^2} F(t,x^i,\Omega,\nu) dV d\Omega d\nu,
\end{align}
represents the number of radiation carriers located at~$x^i$ in the differential volume~$dV$, traveling along the~$\Omega$ direction in the solid angle~$d\Omega$ and that carry energy in frequencies between~$\nu$ and~$\nu + d\nu$. Here~$h$ is Planck's constant and~$c$ is the speed of light. An alternative quantity used to describe radiation is called the specific radiation intensity~$I(t,x^i,\Omega,\nu)$, related to the distribution function as
\begin{align}
	I(t,x^i,\Omega,\nu) = \frac{h^4 \nu^3}{c^2} F(t,x^i,\Omega,\nu).
\end{align}
However, we shall work with the distribution function directly owing to the fact that it is Lorentz-invariant~\cite{MihWei1984}. Throughout the rest of the paper, we will work in units where~$h = c = 1$ and assume the distribution function to be independent of frequency. 

The flux of radiation along the~$x$,~$y$ and~$z$ directions is then defined as
\begin{align}
	F_x = \int_{\mathbb{S}_2} F \, \Omega_x \, d\Omega, && F_y = \int_{\mathbb{S}_2} F \, \Omega_y \, d\Omega, && F_z = \int_{\mathbb{S}_2} F \, \Omega_z \, d\Omega,
\end{align}
where~$\Omega_i$ is the projection of the normalized $3$-momentum~$\vec{\Omega}$ of a carrier, the direction of travel being along the~$i$-th axis. All integrations are performed over the surface of a unit sphere. Another relevant quantity is the energy density of radiation~$E$ defined as
\begin{align}
	E(t,x^i) = \int_{\mathbb{S}_2} F \, d\Omega.
\end{align}
Other integrals with~$\vec{\Omega}$ and~$F$ can be defined, like the radiation pressure tensor, whose components are given by 
\begin{align}
	P_{ij} = \int_{\mathbb{S}_2} F \, \Omega_i \, \Omega_j \, d\Omega, && i = 1,2,3.
\end{align}
The evolution of radiation carriers in special relativity is described by the relativistic Boltzmann equation
\begin{align}
	p^\mu \frac{\p F}{\p x^\mu} = \mathbb{C}[F], && \mu=0,\ldots,3,
\end{align}
where~$p^\mu$ is the $4$-momentum and~$\mathbb{C}$ is a collision term which describes the interaction between radiation and matter. The components of momentum are related to the azimuthal angle~$\phi$ and polar angle~$\theta$ as
\begin{align}
	p^\mu = (1, \cos \phi \sin \theta, \sin \phi \sin \theta, \cos \theta), && \Omega^i = \frac{p^i}{p^0}.
\end{align}
Expressing the right hand side in terms of the emission, scattering and absorption properties of the medium, the Boltzmann equation can be rewritten as
\begin{align} \label{eq:boltzmann}
	\frac{\p F}{\p t} + \Omega^i \frac{\p F}{\p x^i} = \eta - \kappa_a F + \kappa_s \left(\frac{E}{4 \pi} - F \right),
\end{align}
where we have used Einstein's convention for summation over indices. Here~$\eta$ is the radiative emissivity of matter while~$\kappa_a$ and~$\kappa_s$ are the absorption and scattering coefficients, which are related to the inverse of the mean free path. The total extinction coefficient is defined as the sum of absorption and scattering coefficients. The last term on the right hand side of Eq.~\eqref{eq:boltzmann} assumes that any scattering being considered is elastic and isotropic.

\section{The numerical scheme} \label{sec:scheme}
We consider the following ansatz for the distribution function
\begin{align} \label{eq:ansatz}
	F(t,x^i,\Omega) = \sum_{A = 0}^{N-1} F^A(t,x^i) \Psi_A (\Omega) := F^A \Psi_A,
\end{align}
where~$\Psi_A$ are angular basis functions. These are chosen according to the type of scheme. For the~\fpn scheme, they are chosen to be the real spherical harmonics as described in sub-section~\ref{subsub:fpnbasis}.
For~\sn and~\femn these are defined appropriately over a geodesic grid with~$N$ points. This will be explained momentarily.

Substituting this ansatz in the Boltzmann equation in Eq.~\eqref{eq:boltzmann}, we obtain
\begin{align}
	\Psi_A \frac{\p F^A}{\p t} + \Omega^i \Psi_A \frac{\p F^A}{\p x^i} = \mathbb{C}[F].
\end{align}
Multiplying by $\Psi^B = \Psi_B$ and integrating over~$\mathbb{S}_2$ with respect to~$d\Omega$, we obtain
\begin{align}
	M^B_{\ A} \frac{\p F^A}{\p t} + S^{i B}_{\ \ A} \frac{\p F^A}{\p x^i} = \mathbb{S}^{B}[F],
\end{align}
where
\begin{align} \label{eq:massstiff}
	M^B_{\ A} = \int_{\mathbb{S}_2} \Psi^B \Psi_A d\Omega, && S^{i B}_{\ \ A} = \int_{\mathbb{S}_2} \Omega^i \Psi^B \Psi_A d\Omega, && \mathbb{S}^{B}[F] = \int_{\mathbb{S}_2} \Psi^B \mathbb{C}[F] d\Omega, 
\end{align}
are the mass, stiffness and source matrices respectively. Both the mass and stiffness matrices are symmetric matrices independent of position and time and can therefore be pre-computed. The source integral generally has to be evaluated at every time step. Multiplying both sides of Eq.~\eqref{eq:massstiff} by the inverse mass matrix, we obtain a system of coupled advection equations with source terms
\begin{align} \label{eq:boltzmannfinal}
	\frac{\p F^A}{\p t} + \tilde{S}^{i A}_{\ \ B} \frac{\p F^B}{\p x^i} = \tilde{S}^A [F].
\end{align}
For the particular types of sources considered in this paper, parts of the source integral can be pre-computed. The right hand side of Eq.~\eqref{eq:boltzmannfinal} can be written down as
\begin{align}
	\tilde{S}^A[F] = e^A + P^A_{\ B} F^B,
\end{align} 
where the pre-computed parts of the source terms are
\begin{align}
	e^A =  \left[ M^{-1} \right]^{A}_{\ B} \int_{\mathbb{S}_2} \eta \Psi^B d\Omega, && P^A_{\ B} = \frac{1}{4 \pi} \left[ M^{-1} \right] ^{ C}_{ \ B} \int_{\mathbb{S}_2} \Psi_C (\Omega') \ d\Omega ' \int_{\mathbb{S}_2} \kappa_s \Psi^A (\Omega) d\Omega - \left[ M^{-1} \right]^{ C}_{\ B} \delta^{A}_{\ C}.
\end{align}
Here~$\delta^A_{\ C}$ is the Kronecker delta function.
\subsection{Angular discretization} \label{subsec:angular}
An obvious choice for representing discretized angles would be the use of spherical latitude-longitude grids. However, such a choice is not without its disadvantages, the principal of which is the presence of singularities at the poles. This issue can be completely sidestepped by the use of geodesic grids~\cite{Gir1997,HeiRan1995}, which have the additional advantage of being approximately uniformly spaced. While in this paper we make use of uniform grid refinement, geodesic grids have also been used to refine selectively along preferred angular directions~\cite{BohKidTeu2016}.
\begin{figure*}[!t] 
	\centering
	\includegraphics[width=\textwidth]{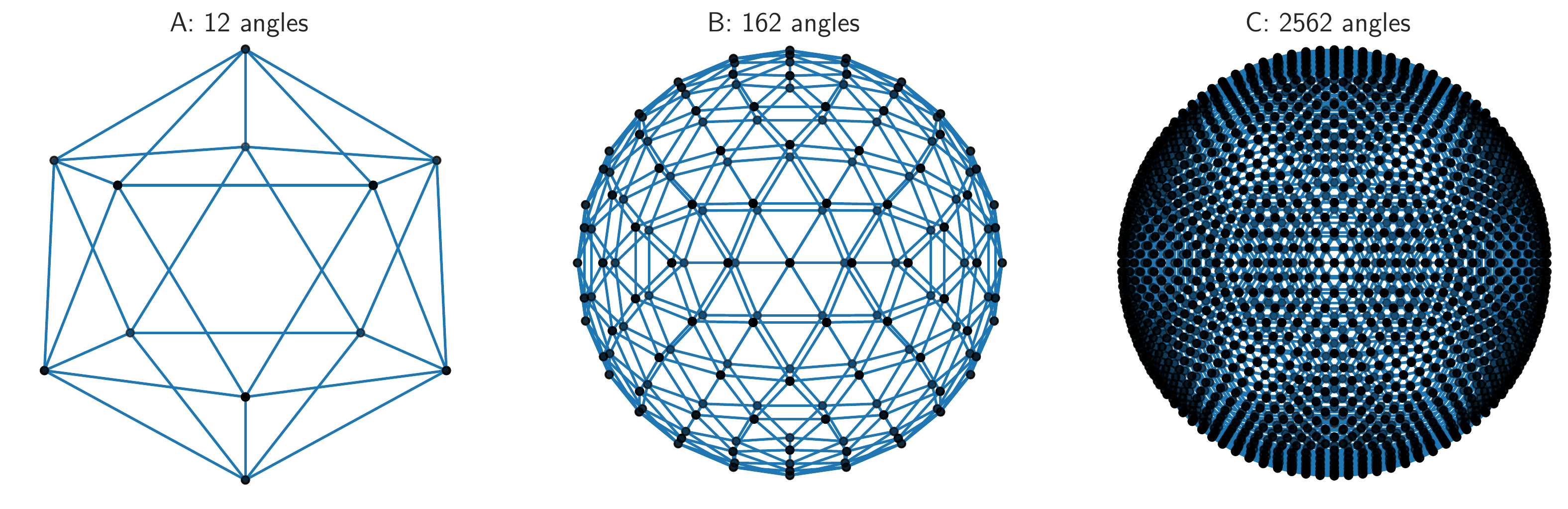}
	\caption{\textit{Left to right}: Plots of the geodesic grid in Cartesian space with varying levels of refinement. Panel A shows the initial icosahedral grid which has~12 vertices,~30 edges, and~20 triangles. Panel B shows the base grid after~2 refinements. This grid has~162 vertices,~480 edges, and~320 triangles. Panel C shows the base grid after~4 refinements. This grid has~2562 vertices,~7680 edges, and~5120 triangles.}
	\label{fig:grid1}
\end{figure*}

The simplest geodesic grid, which we will refer to as the base grid or the unrefined grid, consists of 12 angular points which are vertices of a regular icosahedron, the Cartesian coordinates of which are given by
\begin{align}
	\frac{1}{\sqrt{1+\varphi^2}}(0, \pm 1, \pm \varphi),  && \frac{1}{\sqrt{1+\varphi^2}}(\pm 1, \pm \varphi, 0), && \frac{1}{\sqrt{1+\varphi^2}}(\pm \varphi, 0, \pm 1), && \varphi = \frac{1 + \sqrt{5}}{2}.
\end{align}
Here~$\varphi$ is the golden radio. These vertices lie on the surface of a sphere of unit radius. Each of these vertices have~5 neighboring vertices with which they form a total of~30 unique edges and~20 equilateral triangles~\cite{Wei2022}. Panel A of Fig.~\ref{fig:grid1} shows the initial or base grid.

In the rest of this sub-section, we use Cartesian coordinates~$\vec{x}$ to describe points on the geodesic grid. To refine the spherical grid, we consider an edge of the current grid shared by two grid points with coordinates~$\vec{x}_A$ and~$\vec{x}_B$ and bisect it to obtain coordinates of the midpoint~$\vec{x}_C = (\vec{x}_A + \vec{x}_B)/2$. This point lies inside the unit sphere and its rescaled by a factor of~$1/|\vec{x}_C|$ to project it onto the surface. The process is continued for each edge to obtain the refined grid. Each triangular element of the grid gives rise to~4 smaller triangles. These triangles are no longer equilateral but are almost equilateral and have almost equal areas.

For~$k$ levels of refinement, the number of points, edges and triangles of the geodesic grid are given by~\cite{Gir1997}
\begin{align}
	N_{\textrm{points}} = 12 \times 4^k - 6 \sum_{i=0}^{k-1} 4^i, && N_{\textrm{edges}} = 3 (N_{\textrm{points}}-2), && N_{\textrm{triangles}} = 2 (N_{\textrm{points}}-2), && k \geq 1.
\end{align}
Panels B and C of Fig.~\ref{fig:grid1} show the geodesic after~2 and~4 levels of refinement. They have~642 and~2562 angles respectively.

To compute integrals of functions over these grids, we consider integration over individual triangular elements and then perform a summation over all elements to obtain the desired result. For a single triangular element with three vertices~$A$,~$B$ and~$C$, consider the planar triangle $\Delta_\textrm{p} ABC$ formed by them. Any point~$D$ on or inside this planar triangle can be expressed in terms of three barycentric or areal coordinates~($\xi_1$,~$\xi_2$,~$\xi_3$) defined as~\cite{Wei2022b}
\begin{align} \label{eq:barycentric}
	\xi_1 = \frac{\textrm{area} \ \Delta_\textrm{p} BCD}{\textrm{area} \ \Delta_\textrm{p} ABC}, && \xi_2 = \frac{\textrm{area} \ \Delta_\textrm{p} ACD}{\textrm{area} \ \Delta_\textrm{p} ABC}, && \xi_3 = \frac{\textrm{area} \ \Delta_\textrm{p} ABD}{\textrm{area} \ \Delta_\textrm{p} ABC}.
\end{align}
The three coordinates are not independent but are related by
\begin{align}
	\xi_1 + \xi_2 + \xi_3 = 1.
\end{align}
If the vertices of~$A$,~$B$ and~$C$ in Cartesian coordinates are~$\vec{x}_1$,~$\vec{x}_2$ and~$\vec{x}_3$, then the Cartesian coordinates of any point on the planar triangle with barycentric coordinates~$(\xi_1,\xi_2,\xi_3)$ are given by
\begin{align}
	\vec{x} = \xi_1 \vec{x}_1 + \xi_2 \vec{x}_2 + \xi_3 \vec{x}_3.
\end{align}
The function to be integrated over a triangular element~$\Delta_{\textrm{sp}} ABC$ is first expressed in barycentric coordinates of the planar triangle constructed from the element vertices. The integral over the spherical triangle~$\Delta_{\textrm{sp}} ABC$ is then performed by using the appropriate Jacobian transformation between the planar and spherical triangles.
\begin{figure*}[!t] 
	\centering
	\includegraphics[width=\textwidth]{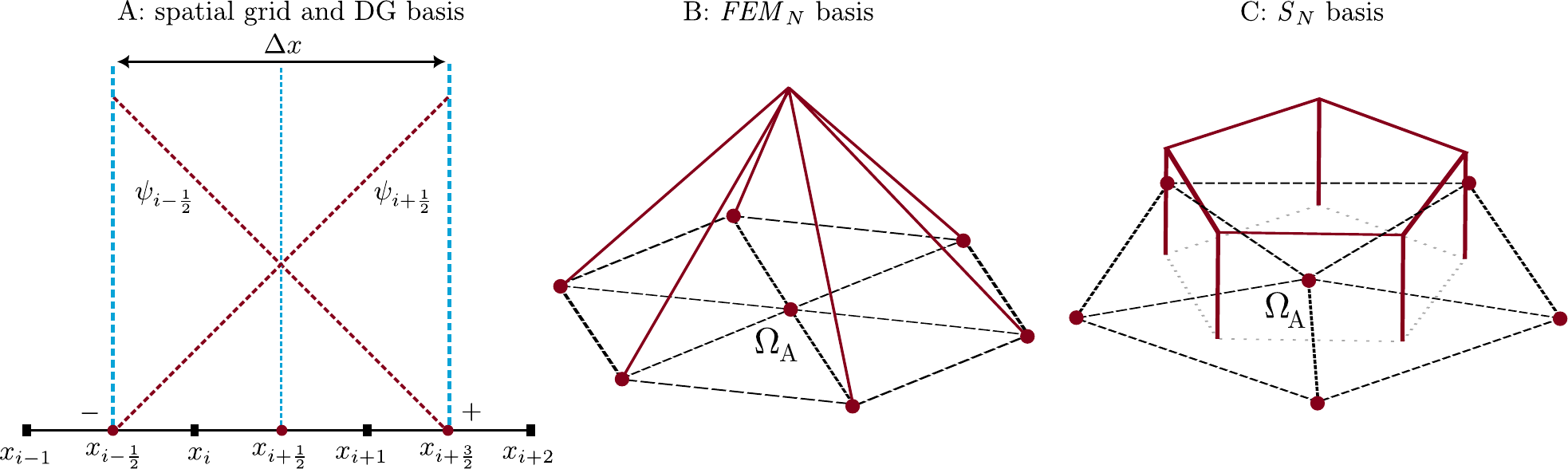}
	\caption{\textit{Left to right}: Basis functions on the spatial and angular grids. Panel A shows the grid structure of a Cartesian element in one dimension. Each element of length~$\Delta x$ is composed of two cells with cell centers located at~$x_i$ and~$x_{i+1}$. The dotted magenta lines are the two Lagrangian basis functions~$\psi_{i-1/2}$ and~$\psi_{i+1/2}$. Panel B shows the~\femn basis function~$\Psi_A$. This function is unity when~$\Omega = \Omega_A$ and vanishes linearly at the neighboring vertices of the grid. Panel C shows a~\sn basis function~$\Psi_A$ which is piece-wise linear and form a non-overlapping pentagon.}
	\label{fig:grid2}
\end{figure*}
\subsubsection{The~\femn scheme} \label{subsub:femnbasis}
In the finite element in angle scheme, we choose finite element basis functions over the geodesic grid~$\Psi_A$ at the vertex~$A$ and which vanishes at the neighboring vertices, as is shown in panel B of Fig.~\ref{fig:grid2}. In a triangular element~$\Delta ABC$ of the grid with barycentric coordinates defined as in Eq.~\ref{eq:barycentric}, the basis function peaking at vertex~$A$ and vanishing at the other two vertices can be written down as
\begin{align}
	\Psi_A(\xi_1,\xi_2,\xi_3) = 2 \xi_1 + \xi_2 + \xi_3 - 1.
\end{align}
In the example plot, the vertex~$A$ has~6 neighbors and therefore six triangular elements must be considered to construct the entire basis at~$\Omega_A$. When vertices~$A$ and~$B$ share an edge, integrals of the product of~$\Psi_A$ and~$\Psi_B$ are non-vanishing, thereby leading to non-zero off-diagonal terms in the mass and stiffness matrices. In practice, we find that our numerical experiments yield better results if we apply mass lumping, that is the original mass matrix in Eq.~\eqref{eq:massstiff} is replaced by a diagonal matrix~$\bar{M}^B_{\ A}$ with elements
\begin{align}
	\bar{M}^C_{\ C} = \sum_{D = 0}^{N-1} M^D_{\ C}.
\end{align}
These basis functions are overlapping and have non-vanishing integrals for off diagonal elements in the mass and stiffness matrices as defined in Eq.~\eqref{eq:massstiff}. While both the lumped mass~$\bar{M}^A_{\ B}$ and stiffness~$S^{i A}_{\ \ B}$ matrices are symmetric, the product of the inverse of the lumped mass with the stiffness matrices~$\tilde{S}^i = [\bar{M}^{-1}] \ S^i$ are not symmetric in general and an eigenvalues of these matrices reveal small imaginary components in them. In practice, the absolute values of the imaginary part of these eigenvalues are~$\sim10^{-14}$ or smaller and are neglected for our calculations.

The eigenspeeds of the system are bounded by the speed of light~$c=1$ and while superluminal modes are absent in the system, instabilities may arise due to the presence of zero speed modes and these must be treated carefully. A prescription for treating zero speed modes is presented in section~\ref{subsection:dg}.
\subsubsection{The~\sn scheme} \label{subsub:snbasis}
An alternative choice of basis results in the~\sn scheme. Inside a similar triangular element~$\Delta ABC$ defined as before, the basis function is
\begin{align}
	\Psi_A(\xi_1,\xi_2,\xi_3) = \begin{cases}
		1,\quad  \xi_1 \geq \xi_2 \ \textrm{and} \ \xi_1 > \xi_3, \\
		0, \quad \textrm{otherwise}.
	\end{cases}
\end{align}
These are non-overlapping ``honeycomb" shaped basis functions as shown in panel C of Fig.~\ref{fig:grid2}. The matrices~$\bar{M}_{AB}$,~$S^i_{AB}$ and~$\tilde{S}^i_{AB}$ are diagonal matrices and $\tilde{S}^i_{AB}$ has eigenvalues bounded by the speed of light. In the case of the~\sn method, the advection operator has the advantage of being diagonal while allows for sweeping in implicit discretization schemes. However, since our problems of interest are relativistic in nature, implicit discretization strategies are not required as advection in such scenarios is non-stiff in nature.
\subsubsection{The~\fpn scheme} \label{subsub:fpnbasis}
An alternative choice to the ansatz in Eq.~\eqref{eq:ansatz} is to choose global basis functions in the form of real spherical harmonics~\cite{RadAbdRez2013}, that is
\begin{align}
	F(t,x^i,\Omega) = \sum_{l,m} F^{lm}(t,x^i) Y_{lm} (\Omega).
\end{align}
These real spherical harmonics can be constructed from the solutions of the Laplace-Beltrami operation on the two sphere~$Y^m_l$
\begin{align}
	Y_{lm} =
	\begin{cases}
		\frac{1}{\sqrt{2}} (Y^m_l + (-1)^m Y^{-m}_l) = \sqrt{2} N^m_l \cos (m \phi) P^m_l (\cos \theta), & m > 0, \\
		Y^0_l, & m = 0, \\
		\frac{1}{i \sqrt{2}} (Y^{-m}_l - (-1)^m Y^{m}_l) = \sqrt{2} N^{|m|}_l \sin (|m| \phi) P^{|m|}_l (\cos \theta), & m < 0,
	\end{cases} 
	\ N^m_l = \sqrt{\frac{(2 l + 1) (l-m)!}{4 \pi (l+m)!}}.
\end{align}
The normalization for~$Y^m_l$ is chosen such that
\begin{align}
	\int_{\mathbb{S}_2} Y_{lm} (\Omega) Y_{l'm'} d\Omega = \delta_{m m'} \delta_{l l'}.
\end{align}
For this basis, the mass matrix becomes the identity, and the stiffness matrix is a symmetric matrix with no superluminal modes. Since the scattering operator in the case of~\fpn is diagonal, this leads to straightforward implicit schemes in scattering dominated regions. However, as we shall see later, this convenience comes at a price, namely, the lack of positivity preservation in certain problems.
\subsection{Spatial discretization} \label{subsection:dg}
We choose an AP~\cite{McCEvaLow2008,LowMor2002} discontinuous Galerkin~(DG) scheme for the spatial discretization. This ensures that the transport equations exhibit the correct behavior when the Knudsen number of the system~\cite{JanLiQiu2014}, that is, the ratio of the mean free path to the characteristic length scale of the problem approaches very small values. In this limit, the Boltzmann equation becomes  diffusion-dominated. To obtain the correct diffusion rates, a non-AP method requires spatial resolution comparable to the mean free path, leading to very high computational cost in these regions. The~AP DG scheme is designed to sidestep this requirement of very high grid resolutions. Solutions from the resulting numerical scheme converge to the true solution of the Boltzmann equation irrespective of the optical depth. Our choice of DG based spatial discretization is described in~\cite{RadAbdRez2013}, from which the main equations are reproduced here for convenience.

We choose a uniformly spaced Cartesian grid with points along each direction being separated by a distance~$(\delta_x, \delta_y, \delta_z)$ such that the coordinates of a point~$(x_i,y_j,z_k)$ is given by $(i \delta_x, j \delta_y, k \delta_z)$. These points represent the centers of a volume cell. An element of the grid in three dimensions is composed of~8 such cells and has dimensional measurements~$(\Delta x, \Delta y, \Delta z) = (2 \delta_x, 2 \delta_y, 2 \delta_z)$. This unique grid structure is chosen to allow for easy integration with finite volume based hydrodynamical codes.

Consider the one-dimensional semi-discrete Boltzmann equation
\begin{align}
	\frac{\p F^A}{\p t} + \tilde{S}^{x A}_{\ \ \ B} \frac{\p F^B}{\p x} = \mathbb{S}^A[F],
\end{align}
The numerical domain is divided into equispaced elements with an element~$\left[x_{i-1/2},x_{i+3/2}\right]$ of the spatial grid shown in panel A of Fig.~\ref{fig:grid2}. The element comprises two cells with cell centers at positions~$x_{i}$ and~$x_{i+1}$ where the distribution function is defined and evolved. Panel A also shows the two linear~DG basis functions~$\psi_{i-1/2}$ and~$\psi_{i+1/2}$ on the element defined by
\begin{align}
	\psi_{i-1/2}(x) = 1 - \frac{x - x_{i-1/2}}{x_{i+3/2} - x_{i-1/2}}, && \psi_{i+3/2}(x) = \frac{x - x_{i-1/2}}{x_{i+3/2} - x_{i-1/2}}.
\end{align}
Given a function~$f$ whose values at the element boundaries~$f_{i-1/2}$,~$f_{i+3/2}$ are known, the basis functions can be used to reconstruct~$f$ anywhere inside the element using
\begin{align}
	f(x) = f_{i-1/2} \psi_{i-1/2}(x) + f_{i+3/2} \psi_{i+3/2}(x).
\end{align}
Values of the evolved variables can be computed at the element edges from their values at the cell centers and vice versa using the relations
\begin{align}
	f_{i-1/2} = \frac{3}{2} f_{i} - \frac{1}{2} f_{i+1}, && f_{i+3/2} = - \frac{1}{2} f_{i} + \frac{3}{2} f_{i+1}.
\end{align}
The numerical scheme then becomes
\begin{align}
	\frac{d F^A_i}{dt} = \frac{1}{\Delta x} \mathbb{F}^A_i, && \mathbb{F}^A_{i} \equiv \frac{3}{2} \mathcal{F}^{-} - \bar{\mathcal{F}} - \frac{1}{2}\mathcal{F}^{+}, && \mathbb{F}^A_{i+1} \equiv \frac{1}{2} \mathcal{F}^{-} + \bar{\mathcal{F}} - \frac{3}{2}\mathcal{F}^{+},
\end{align}
where the flux terms at the cell centers inside an element are functions of an average flux and solutions of the approximate Riemann problem at the two edges. The average flux is defined as
\begin{align}
	\bar{\mathcal{F}} = \frac{1}{2} \left( \tilde{S}^{x A}_{\ \ \ B} F^B_i + \tilde{S}^{x A}_{\ \ \ B} F^B_{i+1}\right),
\end{align}
while the flux at the left edge is given by
\begin{align}
	\mathcal{F}^{-} = \frac{1}{2} \left[ \tilde{S}^{x A}_{\ \ \ B} \left(F^B_L + F^B_R\right) - \hat{S}^{x A}_{\ \ \ B} \left(F^B_{R} - F^B_{L}\right)\right].
\end{align}
Here $F^B_L$ and $F^B_L$ are the left and right states of the Riemann problem at the left edge. The matrix $\hat{S}^{x A}_{\ \ \ B}$ in the second term controls the amount of numerical dissipation introduced to effectively tackle zero speed mode instabilities and is defined as
\begin{align}
	\hat{S}^{x A}_{\ \ \ B} = \mathcal{R}^{x A}_{\ \ \ C} \max(v,|\Lambda^{x C}_{\ \ \ D}|) \mathcal{L}^{x D}_{\ \ \ B},
\end{align}
where~$\mathcal{L}$ and $\mathcal{R}$ are the matrices formed from the set of left and right eigenvectors of~$\tilde{S}^x$ and~$\Lambda^x$ is a diagonal matrix of their corresponding eigenvalues. $v$ is a positive dissipation parameter which introduced artificial dissipation, which, in our numerical experiments was chosen to be~$1/\sqrt{3}$. In terms of the cell centered values in an element and its neighbors,~$\mathcal{F}^{-}$ and~$\mathcal{F}^{+}$ can be written as
\begin{align}
	\mathcal{F}^{-} &= \frac{1}{2} \left[ \tilde{S}^{x A}_{\ \ \ B} \left(-\frac{1}{2}F^B_{i-2} + \frac{3}{2} F^B_{i-1} + \frac{3}{2} F^B_{i} - \frac{1}{2} F^B_{i+1}\right) - \hat{S}^{x A}_{\ \ \ B} \left(-\frac{1}{2}F^B_{i-2} + \frac{3}{2} F^B_{i-1} - \frac{3}{2} F^B_{i} + \frac{1}{2} F^B_{i+1}\right)\right], \nonumber \\
	\mathcal{F}^{+} &= \frac{1}{2} \left[ \tilde{S}^{x A}_{\ \ \ B} \left(\frac{3}{2}F^B_{i+2} - \frac{1}{2} F^B_{i+3} - \frac{1}{2} F^B_{i} + \frac{3}{2} F^B_{i+1}\right) - \hat{S}^{x A}_{\ \ \ B} \left(\frac{3}{2}F^B_{i+2} - \frac{1}{2} F^B_{i+3} + \frac{1}{2} F^B_{i} - \frac{3}{2} F^B_{i+1}\right)\right].
\end{align}

In order to prevent the appearance of artificial extrema in numerical solutions, we employ the use of slope limiters. We consider an element comprising the~$i$ and~$i+1$ cells, At every sub-step of the time integrator, the average value of the distribution function in this element~$\bar{F}^A_{i,i+1}$ and the values at the cell centers~$F^A_i$,~$F^A_{i+1}$ are used to perform a reconstruction procedure to obtain new cell-centered values~$\tilde{F}^A$,
\begin{align}
	\tilde{F}^A_i = \bar{F}^A_{i,i+1} + \sigma_x (x_i - \bar{x}), && \tilde{F}^A_{i+1} = \bar{F}^A_{i,i+1} + \sigma_x (x_{i+1} - \bar{x}),
\end{align}
where~$\bar{x}$ is the mean coordinate position of the element and~$\sigma_x$ is the slope limited correction from a limiter function defined as
\begin{align}
	\sigma_x = \mathcal{M}\left(\frac{F^A_{i+1} - F^A_{i}}{\Delta x /2},  \frac{\bar{F}^A_{i,i+1}-\bar{F}^A_{i-2,i-1}}{\Delta x},  \frac{\bar{F}^A_{i+2,i+3}-\bar{F}^A_{i,i+1}}{\Delta x} \right),
\end{align}
where~$\mathcal{M}$ is chosen to be the minmod limiter defined as~\cite{Lev2002}
\begin{align}
	\textrm{minmod}(a,b,c) = \begin{cases}
		\textrm{min}(|a|,|b|,|c|), & \textrm{sign}(a) = \textrm{sign}(b) = \textrm{sign}(c), \\
		0, & \textrm{otherwise},
	\end{cases}
\end{align}
or the ``sawtooth free" double minmod limiter~\cite{McCLow2008}
\begin{align}
	\textrm{s-minmod2}(a,b,c) = \begin{cases}
		s \times |a|, & |a| < 2 \times \textrm{min}(|b|,|c|), \\
		s \times \textrm{min}(|b|,|c|), & \textrm{otherwise},
	\end{cases}, \ 
	s = \begin{cases}
		\textrm{sign}(a), & \textrm{sign}(a) = \textrm{sign}(b) = \textrm{sign}(c), \\
		0, &\textrm{otherwise}. \\
	\end{cases} 
\end{align}
The s-minmod2 limiter is asymptotic preserving and does not affect the quality of the solution away from the location of the  extrema~\cite{McCLow2008}. In certain situations, a modified version of the s-minmod2 limiter proves to be more effective in removing negative values. This is defined as
\begin{align}
	\textrm{modminmod2}(a,b,c) = \textrm{s-minmod2}(a,b/2,c/2).
\end{align}
\subsection{Time integration} \label{subsec:time}
For the time integration, we use a second-order Runge-Kutta method. Consider the form of the Boltzmann equation in Eq.~\eqref{eq:boltzmannfinal}
\begin{align}
	\frac{\p F^A}{\p t} = - \tilde{S}^i_{BA} \frac{\p F^A}{\p x^i} + e^A + P^A_{\ B} F^B.
\end{align}
Given a solution~$F^A_k$ of the above equation at time $t = k \Delta t$, we compute the solution at time $t + \Delta t = (k+1) \Delta t$ as
\begin{align}
	F^A_{k+1/2} &= F^A_{k }- \frac{\Delta t}{2} \left( \tilde{S}^i_{BA} \frac{\p F^A}{\p x^i}\Bigg|_k - e^A_k - P^A_{\ B} F^B_{k} \right), \\
	F^A_{k+1} &= F^A_{k} -\Delta t \left(\tilde{S}^i_{BA} \frac{\p F^A}{\p x^i}\Bigg|_{k+1/2} -  e^A_{k+1/2} - P^A_{\ B} F^B_{k+1/2} \right).
\end{align}
An alternative semi-implicit time stepping algorithm as described in~\cite{McCEvaLow2008} may be used for the optically thick regime.

\subsection{Positivity preserving strategies} \label{subsec:limfilt}
The distribution function of carriers or their energy in a system cannot be negative. Different discretizations of the transport equations may not respect this fact, which is why corrections to numerical solutions have to be made to ensure non-negativity. For the~\femn scheme, we use a simple limiting strategy and for~\fpn schemes, a filtering strategy as described in~\cite{RadAbdRez2013} is used.
\subsubsection{Limting} \label{subsub:limting}
Several limiters have been proposed in~\cite{LaiHau2019} for the~\fpn scheme to ensure non-negativity, all of which can be repurposed to act as limiters for the~\femn scheme. For this paper, we choose the clipping limiter~(clp) as described in~\cite{ShaWen2004,LigDur2016}. This truncates the negative values of~$F^A_i$ to zero and compensates for it by readjusting other positive values of~$F^A_i$ by rescaling them by a parameter~$\theta$ defined as
\begin{align}
	\theta = \frac{\sum_{A,B} M_{AB} F^A_{i}}{\sum_{A,B} M_{AB} \tilde{F}^A_{i}}, && \tilde{F}^A_{i} = \textrm{max}(F^A_{i},0).
\end{align}
The values of~$F^A_i$ after limiting becomes
\begin{align}
	 F^{A(\textrm{new})}_i = \begin{cases}
		\theta F^A_{i}, & \textrm{if} \ F^A_{i} > 0, \\
		0, & \textrm{if} \ F^A_{i} \leq 0.
	\end{cases}
\end{align}
An advantage of this limiter is that it conserves~$E$ point-wise.

\subsubsection{Filtering}
For the~\fpn method, we use filtering to remove the effects of Gibbs' oscillations which can drive the solution to acquire non-physical values. At each sub-step of the time integration, a filtering operation is performed on~$F$ as
\begin{align} \label{eq:opacityeqn}
	F^{\textrm{new}} = \sum_{A} \sigma\left( \frac{l}{l_{\textrm{max}}+1}\right)^s F^A Y_A, && s = - \frac{\Delta t \ \sigma_{\textrm{eff}}}{\log \sigma(l/(l+1))},
\end{align}
where the summation operation is performed over $A = \{l,m\}$ and~$\sigma$ is the filter chosen which has a strength~$s$. For most of our numerical experiments, instead of directly specifying the strength, we specify another parameter~$\sigma_{\textrm{eff}}$, called the effective opacity which is related to filter strength as described in Eq.~\ref{eq:opacityeqn}. For all our tests, we use the Lanczos filter
\begin{align}
	\sigma_{\textrm{L}}(x) = \frac{\sin x}{x}.
\end{align}
The filtering strategy is described extensively in~\cite{RadAbdRez2013}.
\begin{figure*}[!th] 
	\centering
	\includegraphics[width=0.97\textwidth]{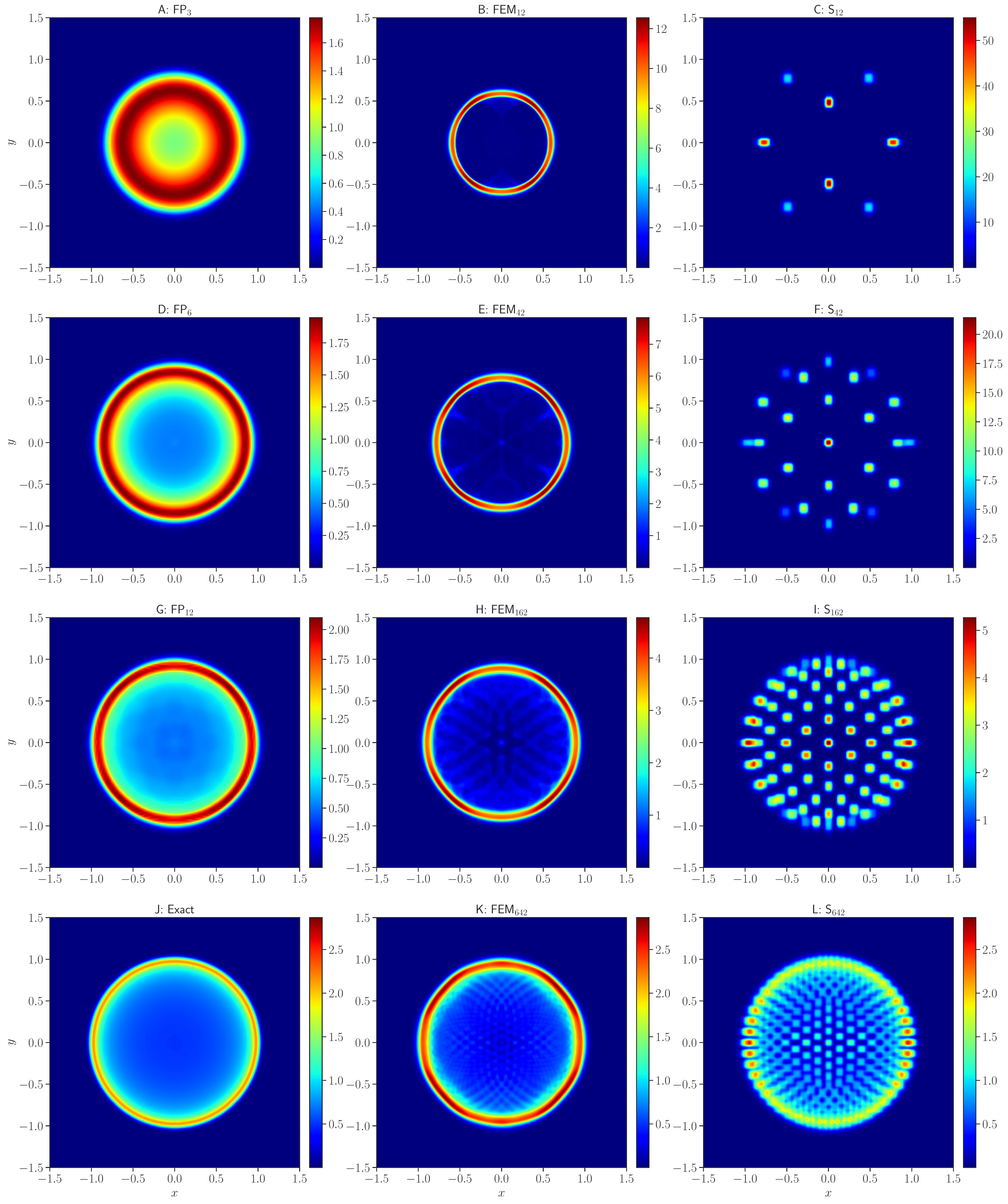}
	\caption{A comparison of~$E$ between the~FEM${}_{N}$,~\sn and~\fpn methods at different angular resolutions. The angular grid for~\femn and~\sn runs have~$12$,~$42$,~$162$ and~$642$ points. For~\fpn runs, we consider~$N = 3$, $6$ and $12$. All simulations are evolved till $t = 1$.  Both $\textrm{\textit{FEM}}_N$ and $S_N$ results approximate the exact solution better with increasing angular resolution, but~\femn is superior for the same number of angles. The~\fpn results approximate the exact solution better than the other methods at intermediate resolutions provided that the filter parameter is chosen appropriately.}
	\label{fig:line1}
\end{figure*}
\section{Numerical results} \label{section:results}
\subsection{The line source test}
\begin{figure*}[!t] 
	\centering
	\includegraphics[width=\textwidth]{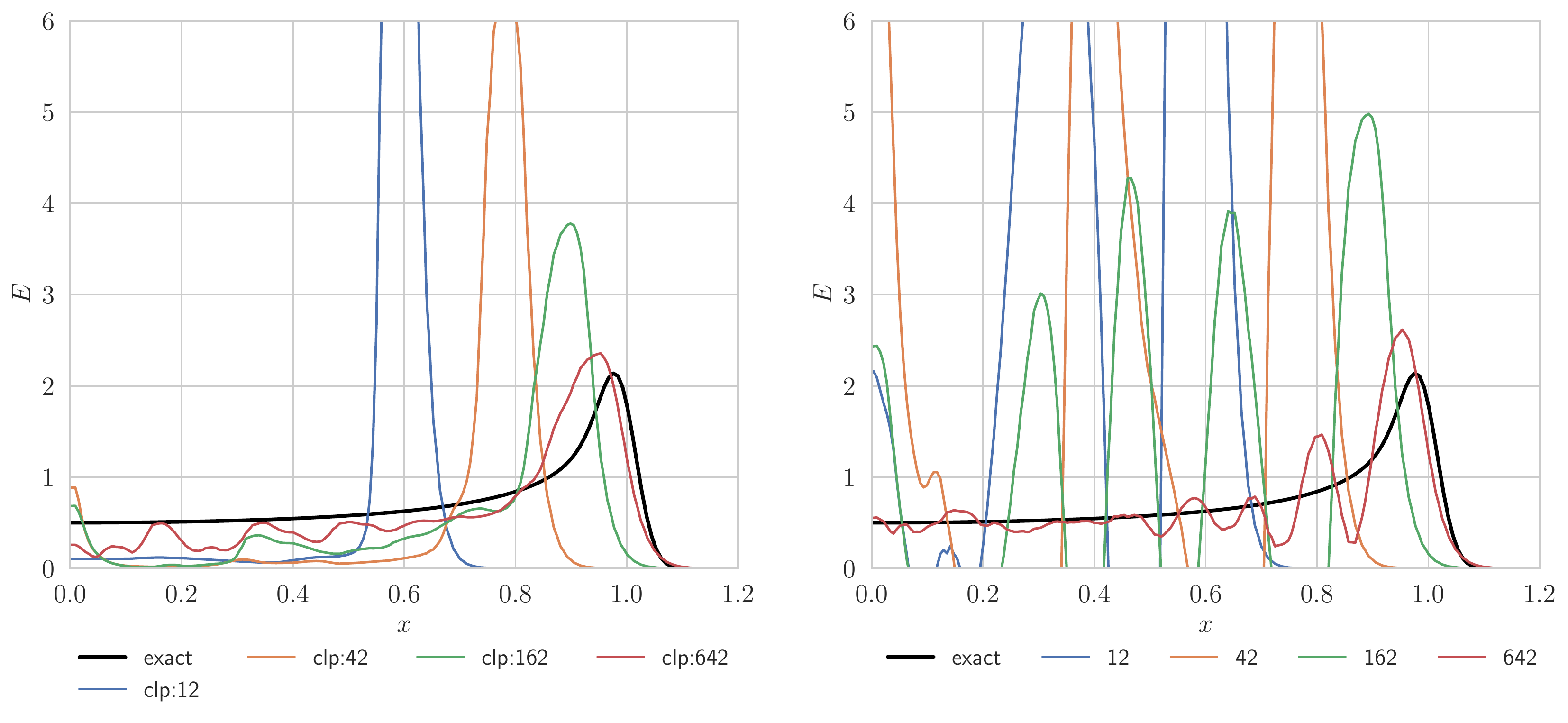}
	\caption{\textit{Left}: A comparison of~$E$ between the exact solution and the numerical solutions along the x-axis with the $\textrm{\textit{FEM}}_N$ method with limiting at 4 angular resolutions. Solutions with the clipping limiter reproduce the exact solution more faithfully in the~$L^1$/$L^{\infty}$ norm with increasing number of angles. \textit{Right}: A comparison of~$E$ along the x-axis between the exact and numerical solutions computed with the~\femn scheme without limiting. For small angles, large oscillations are seen and both $E$ and $F^A$ can acquire negative values. With sufficient angular resolution, the solution becomes non-negative even without limiting.}
	\label{fig:line2}
\end{figure*}
The line source test, proposed in~\cite{Gan1999}, is an important numerical experiment to demonstrate the strengths and weaknesses of radiation transport schemes. We perform this test with our aim being to test the effectiveness of the new angular discretization scheme. The problem consists of a pulse of radiation concentrated along the z-axis which propagates isotropically in vacuum. The initial data for the problem is given by~\cite{GarHau2013}
\begin{align}
	\tilde{F}(t,x,y,\Omega) = \frac{1}{4 \pi} \delta(x, y).
\end{align}
This problem can be solved analytically. The time evolution of the radiation energy density~$E$ is
\begin{align}
	\tilde{E}(t,x,y) = \frac{1}{2 \pi} \frac{H(t-r)}{t \sqrt{t^2 - r^2}}, && r = \sqrt{x^2+y^2},
\end{align}
where $H$ is the Heaviside step function and~$r$ is the cylindrical radius coordinate. The solution can be described as a ``cylindrical shell" of radiation propagating at the speed of light, while in its interior, radiation falls as a function of $r$. In our numerical implementation, we use Cartesian coordinates for the spatial domain, to compare the rotational invariance, or lack thereof, of the~\sn,~\femn and~\fpn schemes.

The delta function in our numerical implementation is represented by a sharp Gaussian centered at the origin whose steepness is controlled by a parameter~$\omega$ chosen according to the prescription in~\cite{GarHau2013}
\begin{align}
	F(t=0,x,y,\Omega) = \textrm{max} \left( \frac{1}{8 \pi \omega^2} e^{-(x^2+y^2)/(2 \omega^2)}, 10^{-4}\right), && \omega = 0.03.
\end{align}
The choice of this type of initial data ensures that artifacts which arise during time evolution are due to the angular discretization of our scheme and not of the spatial discretization. The exact time evolution of $E$ can then be computed by performing a convolution over the initial data
\begin{align} \label{eq:lineexact}
	E(t,x,y) = \int_{\mathbb{R}^2} E(t=0,x,y) \tilde{E}(t,x-x',y-y') dx' dy'.
\end{align}
In our numerical experiments, we choose a spatial grid $x \in \left[-1.5,1.5\right]$, $y \in \left[-1.5,1.5\right]$ with $\delta x = \delta y \approx 0.006$. All simulations are evolved till $t = 1$ with $\delta t = 0.002$. For the~\femn and~\sn schemes, we choose~$12$,~$42$,~$162$ and~$642$ angles for the angular grid and consider both with and without limiting cases for~\femn. For the~\fpn scheme, we consider the values of $N$ to be~$3$,~$6$ and~$12$, which correspond to~$16$,~$49$ and~$169$~$(l,m)$ modes. The Lanczos filter used for~\fpn has $\sigma_{\textrm{eff}}$ as~$20$, following the suggestion in~\cite{RadAbdRez2013}.
\begin{figure*}[!t] 
	\centering
	\includegraphics[width=\textwidth]{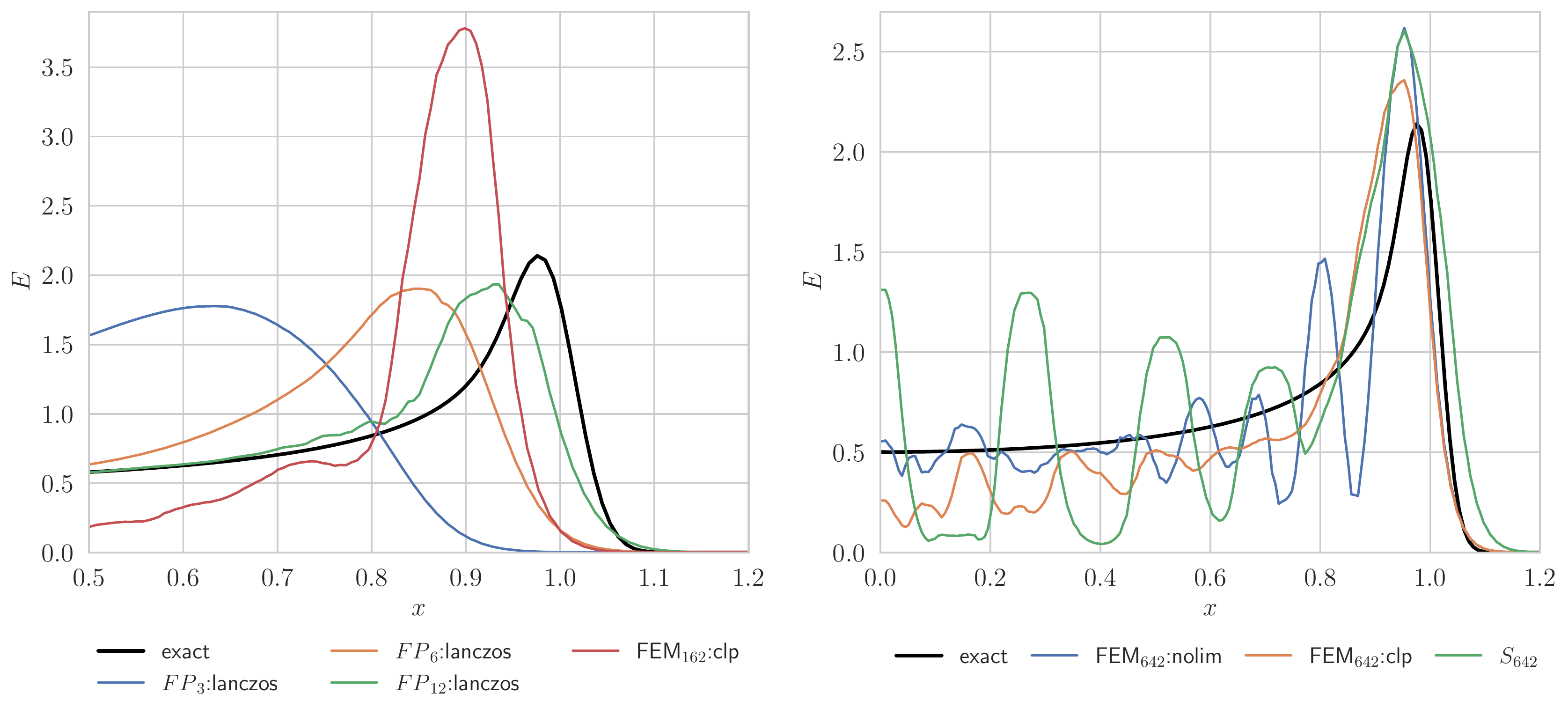}
	\caption{\textit{Left}: A comparison of $E$ between the $\textrm{\textit{FP}}_{N}$ and $\textrm{\textit{FEM}}_{N}$ solutions along the x-axis. The~$FP_{12}$ solution with~$169$~$\left(l,m\right)$ modes is superior to the~FEM${}_{162}$ solution, which is of comparable resolution. \textit{Right}: A comparison of~$E$ between the non-limited, clipped and $S_N$ solutions along the x-axis for 642 angles. The non-limited~\femn solution at this resolution does not contain negative values.}
	\label{fig:line3}
\end{figure*}
\begin{figure*}[!t] 
	\centering
	\includegraphics[width=\textwidth]{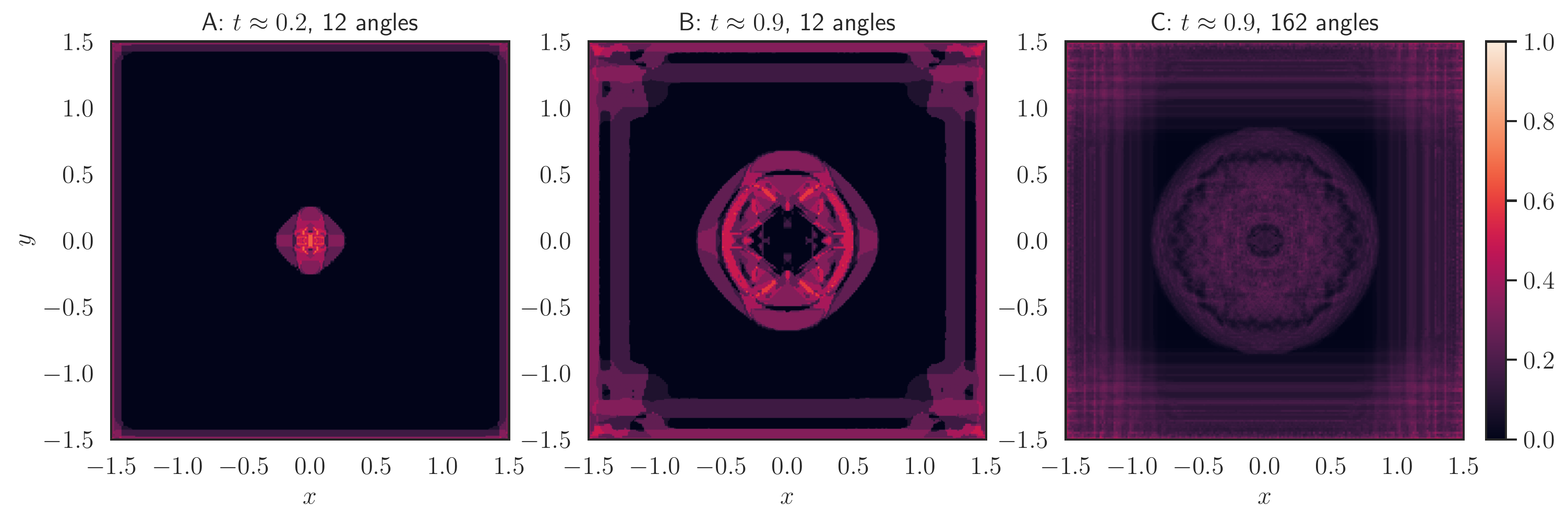}
	\caption{A comparison of the limiter indicator for three different test cases of the line test with the clipping limiter. \textit{Panel A}: The limiter indicator early on during the simulation for the FEM${}_{12}$ case at~$t \approx 0.2$.  \textit{Panel B}: The indicator for the same simulation at~$t \approx 0.9$. \textit{Panel C}: The indicator at~$t \approx 0.9$ for a higher resolution run with~$162$ angles.}
	\label{fig:linelimiter}
\end{figure*}

Without the use of a limiter, the~\femn solutions exhibits oscillations and~$E$ may contain non-physical values at early times, although these are reduced at higher angular resolutions, as can be seen from panel B of Fig.~\ref{fig:line2}. The exact solution for~$E$ resembles a ``cylindrical shell" propagating with a speed of~$1$. The same figure shows that at very low number of angles, this "shell" propagates slower than the speed of light but the numerical solution converges to the true solution of the problem with increase in~$N$. For~$N = 642$, the radiation energy density at~$t = 1$ reproduces the true solution without any negative values. The use of the limiter resolves the problem of negative regions in the solution irrespective of the number of angles considered. This is seen in panel A of Fig.~\ref{fig:line2}. Panels B, E, H and K of Fig.~\ref{fig:line1} show that solution with the clipping limiter retain the ``cylindrical shell" feature even at very low angular resolutions. When compared with the exact solution in panel J of the same figure, we see that the propagation speed of the ``shell" become progressively more accurate with increasing angular resolution. For very low number of angles, the limiter actively removes negative values in the interior of the ``shell", but limiting is required less frequently for higher~$N$. Fig.~\ref{fig:linelimiter} plots an indicator function to show regions of the numerical domain where solutions have become non-physical. This indicator counts the number of points in the phase space which have developed non-physical values and normalizes the result by dividing by the total number of points. A higher number implies that the limiter has to work more aggressively to `correct' for the distribution function.  Panels A and B are simulations of the line test with 12 angles at two different times,~$t \approx 0.2$ and~$t \approx 0.9$. Panel C evolves the line test till~$t \approx 0.9$ but with 162 angles. A comparison between panels B and C demonstrates that less work is needed by the limiter to fix non-physical values in the solution when angular resolution is increased. The line test is an extreme case to test numerical schemes for non-physical oscillations, but in a more realistic scenario, we expect the limiter to be used less frequently.

The~\sn solutions do not exhibit negative values but lack the ``cylindrical shell" feature of the exact solution. Instead, radiation propagates outward as localized blobs at various speeds, bounded by the speed of light. This presence of ``ray effects" can be diminished by increasing the number of angles, as can be seen in panels C, F, I and L of Fig.~\ref{fig:line1}. This makes the~\sn solutions inferior to their~\femn and~\fpn counterparts both in terms of reproducing key features and maintaining rotational invariance.
\begin{figure*}[!t] 
	\centering
	\includegraphics[width=0.99\textwidth]{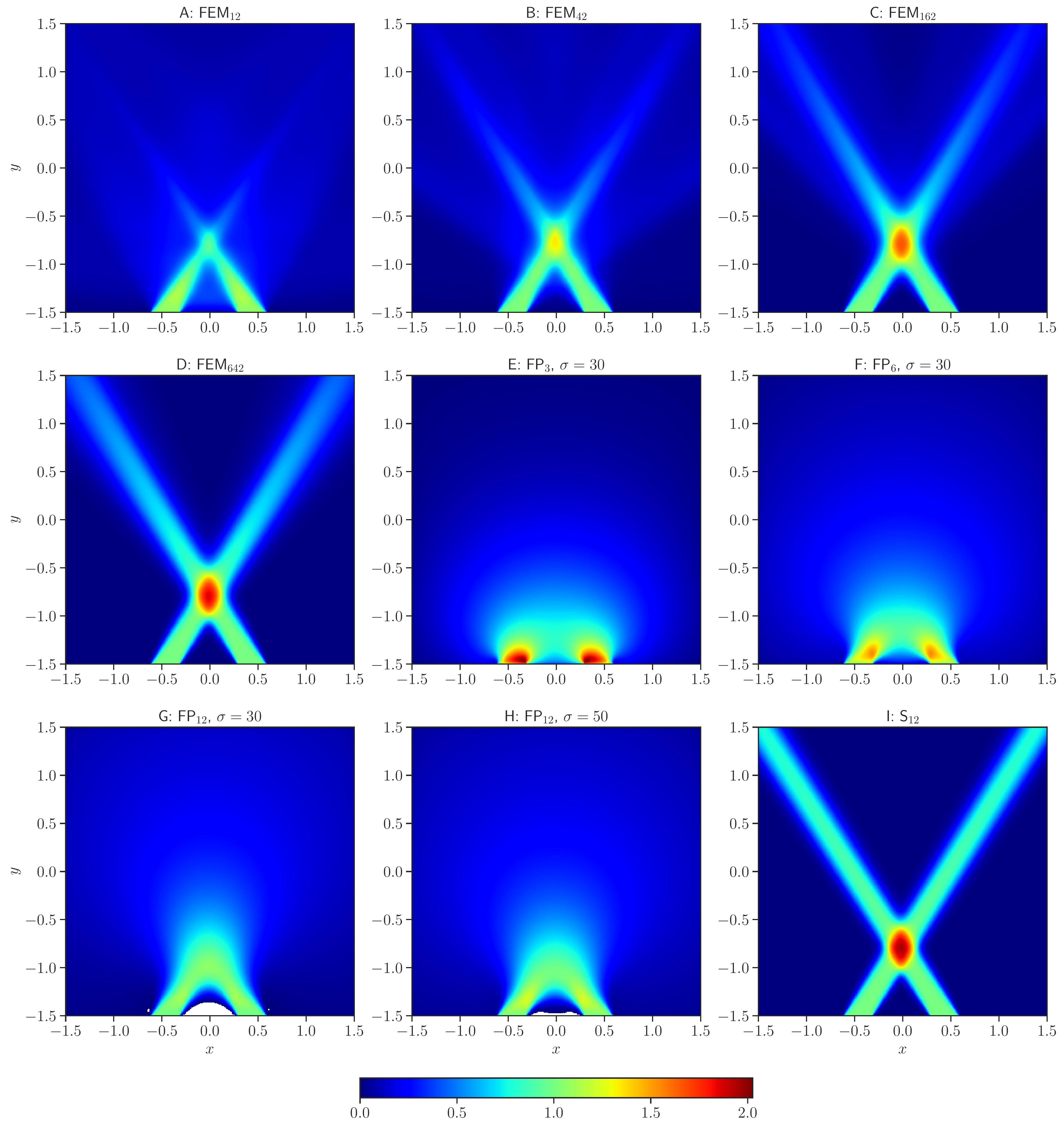}
	\caption{Plots of the searchlight test where two beams enter from the lower boundary of the numerical domain and intersect with each other. These simulations have been evolved up to steady state. The solution to this problem is most accurately represented by the~\sn solution even at very low angular resolution provided the direction of propagation of the beam is a point on the discrete grid. Panels A-D show the~\femn solutions with the limiter at different angular resolutions. Although these solutions become progressively better with the number of angles, they are inferior in quality to their~\sn counterparts. The~\fpn runs develop negative values in the solution irrespective of the choice of~$N$ unless a very high filter opacity is chosen. This in turn produces results inferior to both~\femn and~$S_N$, as can be seen in panels E-H. For~$N = 12$, negative values develop in the solution even when~$\sigma_{\textrm{eff}} = 50$. The white regions in panels G and H are the regions where~$E$ develops non-physical values. A comparison with results from the~$S_{12}$ scheme is provided in panel I.}
	\label{fig:searchlight1}
\end{figure*}

The~\fpn solutions exhibit rotational invariance owing to the properties of spherical harmonics. However spurious oscillations in the solutions have to be remove with the help of a filter with a tunable parameter for its strength. With our choice of filter opacity, the solutions do not exhibit negative values. The low resolution~$\textrm{\textit{FP}}_3$ solution is not able to reproduce the correct speed of propagation of the outer wavefront. The~$\textrm{\textit{FP}}_6$ solution of the other hand is more accurate in representing the exact solution than~\femn solutions possessing approximately the same number of angles. In both~$L^1$ and~$L^{\infty}$ norms, the~$\textrm{\textit{FP}}_{12}$ solutions are superior to~\femn solutions, although at large number of modes or angles, both methods can represent the exact solution with reasonable accuracy. The~\fpn solutions with small number of modes develop negative values unless a filter with a high opacity parameter is chosen. This motivates the choice of the~\femn scheme at low resolutions. However at intermediate number of~$\left(l,m\right)$ modes, the~\fpn solutions are of a superior quality both in terms of accuracy and strict rotational invariance, and therefore are preferred.
\subsection{The searchlight beam test}
Next, we focus our attention to the searchlight beam test~\cite{StoMihNor1992,SumYam2012,PerPenNov2014} in two dimensions. The test involves a narrow beam of radiation propagating in vacuum in a particular direction from the boundary of the numerical domain. In our case, we consider two sources of radiation emitting beams of radiation from the~$y = - 1.5$ boundary at polar angles~$\phi_1 \approx 58.28\degree$ and~$\phi_2 \approx 121.72 \degree$, the angles so chosen that the two beams meet each other inside the numerical domain after a certain period of time. The numerical domain extends from~$-1.5$ to~$1.5$ along each dimension with~$\delta x = \delta y = 0.0075$ and the system is evolved up till~$t = 10$ with~$\delta t \approx 0.0067$ when steady state has been achieved. In this test, we choose the modified s-minmod2 limiter to prevent the appearance of negative values in the solution. The beams should ideally propagate without dispersion and cross each other without interaction. This test is challenging for radiation transport codes due to the presence of sharp gradients in the solution that can give rise to negative values in~$F$ and~$E$.

Fig.~\ref{fig:searchlight1} compares the results of this test between the~$\textrm{FEM}_N$,~\sn and~\fpn schemes. We find that the~\fpn methods are inefficient in depicting beams of radiation without developing non-physical values in the solution. This issue can be resolved by choosing an extremely high opacity parameter for the filter, which in turns affects the quality of the solution adversely. This is clearly seen in panels E and F of Fig.~\ref{fig:searchlight1} where~$\sigma_{\textrm{eff}}$ is chosen to be~$30$ to remove non-physical values in~$E$. For the~$FP_{12}$ case, panels G and H show that even with higher opacities,~$E$ may develop negative values. The~\fpn scheme performs significantly worse when compared with the other two schemes. 

Plots A-D show results with the~\femn with increasing angular refinement, demonstrating that the results get progressively better with increase in~$N$. In all cases, this scheme performed better than the~\fpn scheme. However, the~\sn scheme produces more accurate results than either~\femn or~\fpn irrespective of the level of refinement when the direction of beam propagation lie on our chosen geodesic grid.
\begin{figure*}[!t] 
	\centering
	\includegraphics[width=0.4\textwidth]{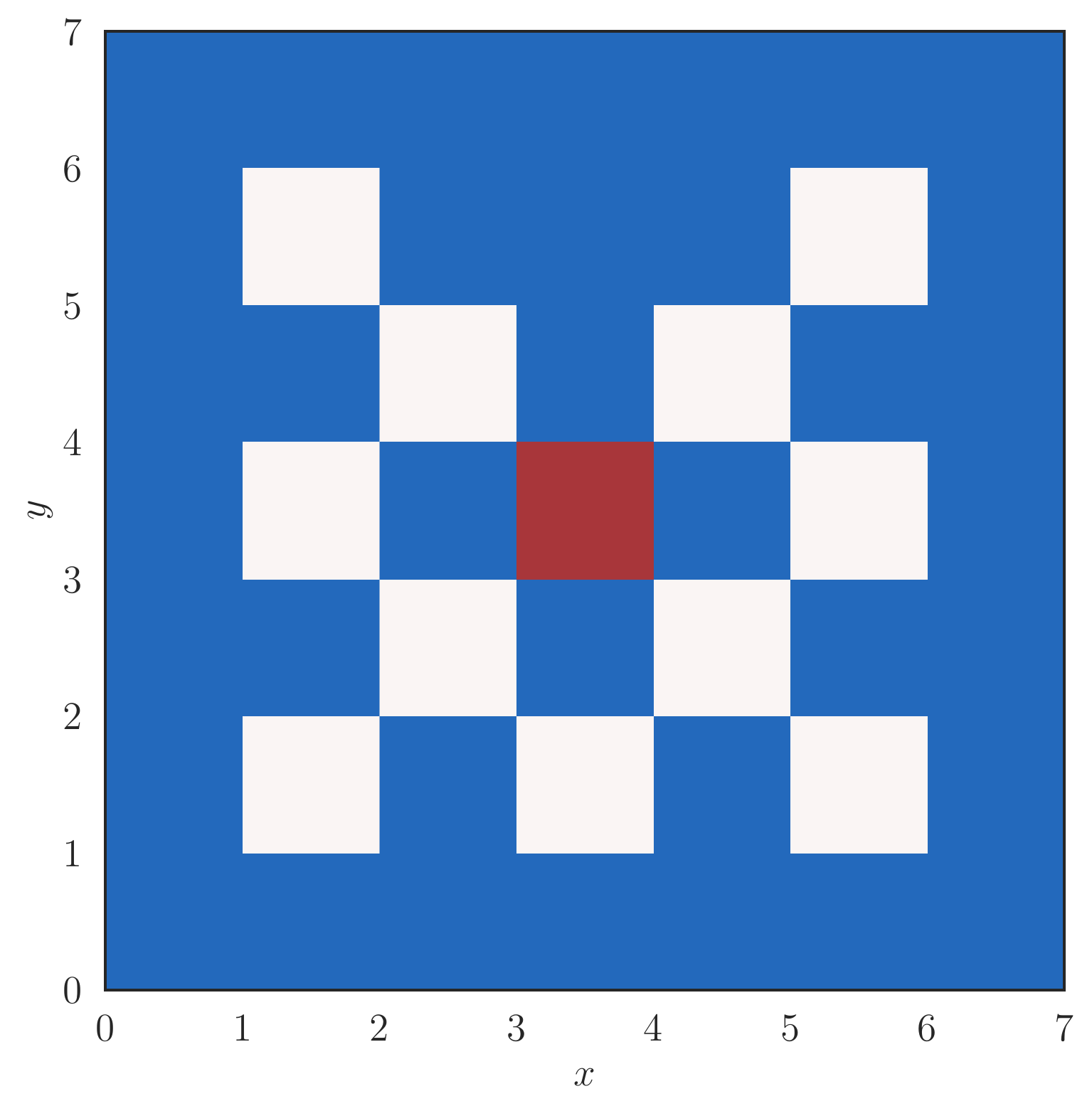}
	\caption{The initial setup of the checkerboard or the lattice problem. The central red square is emitting with~$\eta = 1/4\pi$ and has a scattering opacity~$\kappa_s = 10$. The white squares are absorption regions with~$\kappa_a = 1$, while the surrounding blue area is purely scattering with~$\kappa_s = 1$.}
	\label{fig:lattice0}
\end{figure*}
\begin{figure*}[!t] 
	\centering
	\includegraphics[width=\textwidth]{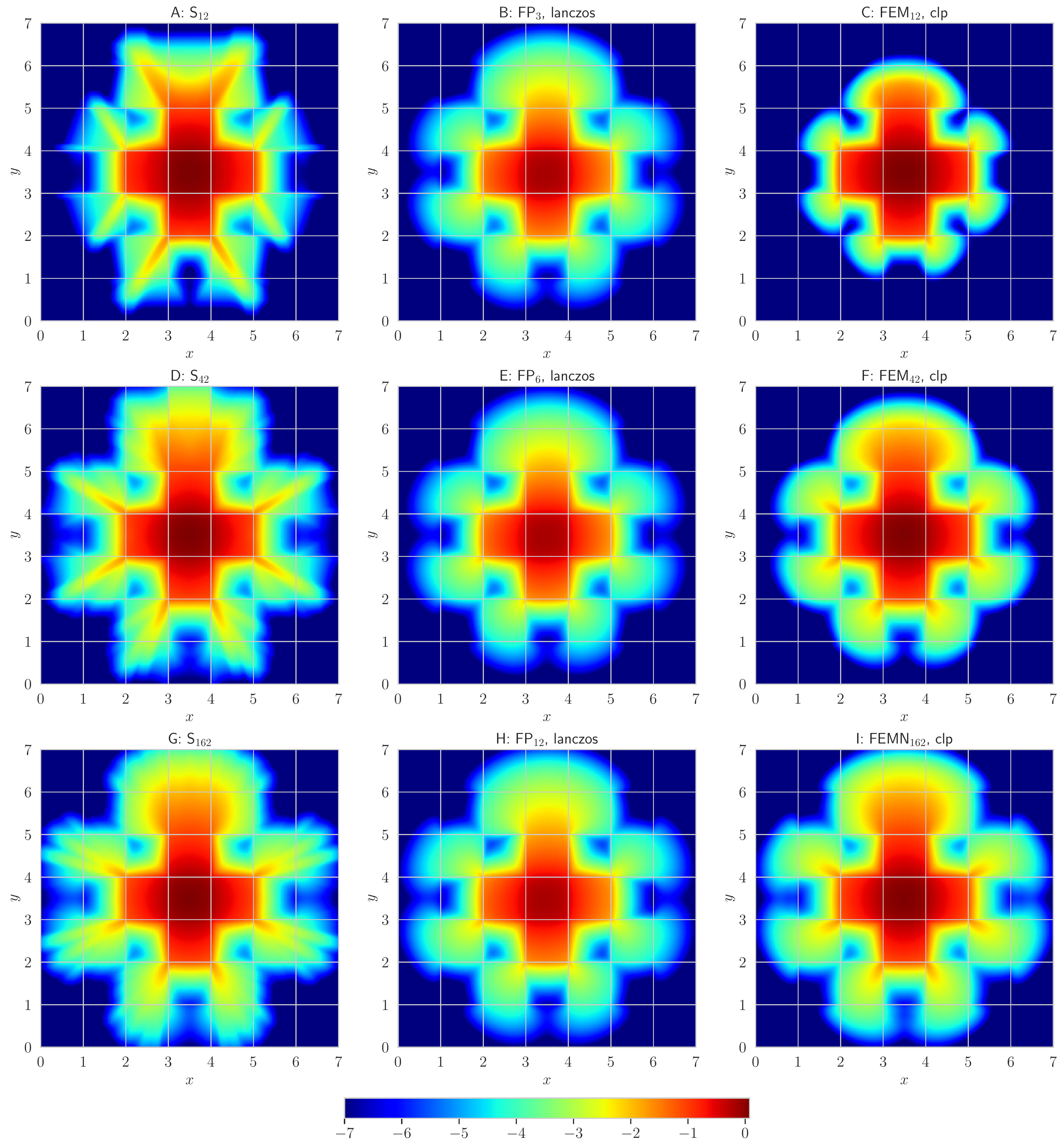}
	\caption{A comparison of $\log_{10} E$ for the lattice test problem evolved up to $t=3.2$ between the~$\textrm{\textit{S}}_{N}$,~\fpn and~\femn methods. At this time, radiation reaches the edges of the outer boundary of the numerical domain. The~\fpn and~\femn schemes produce comparable results with solutions propagating at the correct speeds at moderate angular resolutions. The~\sn method produces prominent ray artifacts at the angular resolutions considered.}
	\label{fig:lattice1}
\end{figure*}
\begin{figure*}[!t] 
	\centering
	\includegraphics[width=\textwidth]{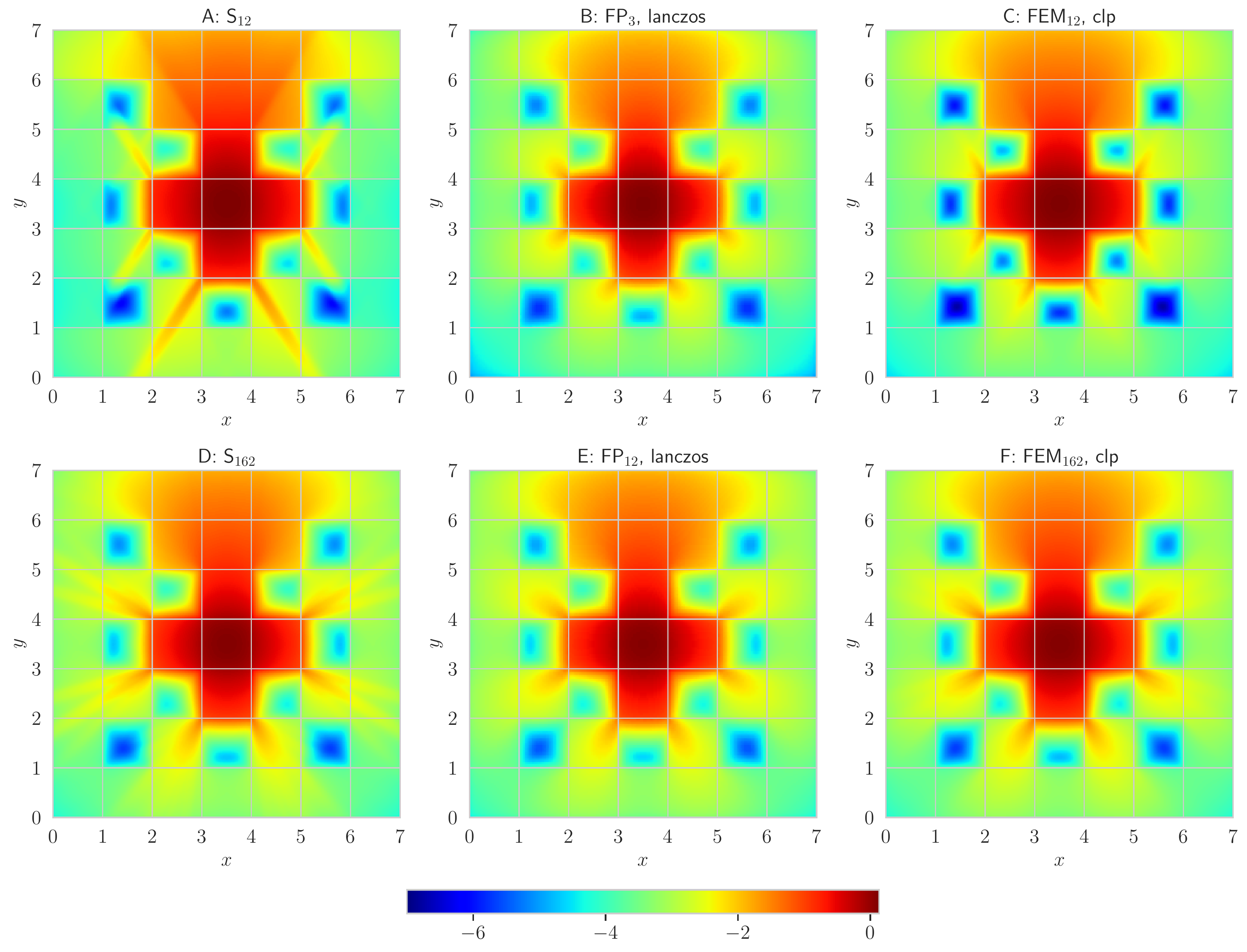}
	\caption{A comparison of results from the three difference schemes for the searchlight problem which has been evolved up to steady state. The solutions for the~\femn and~\fpn schemes are nearly identical, irrespective of angular resolution. The~\sn results still shows significant ray effects at low number of angles. These effects become less prominent with the additional of more angular points.}
	\label{fig:lattice2}
\end{figure*}
\subsection{The lattice test}
\begin{figure*}[!th] 
	\centering
	\includegraphics[width=\textwidth]{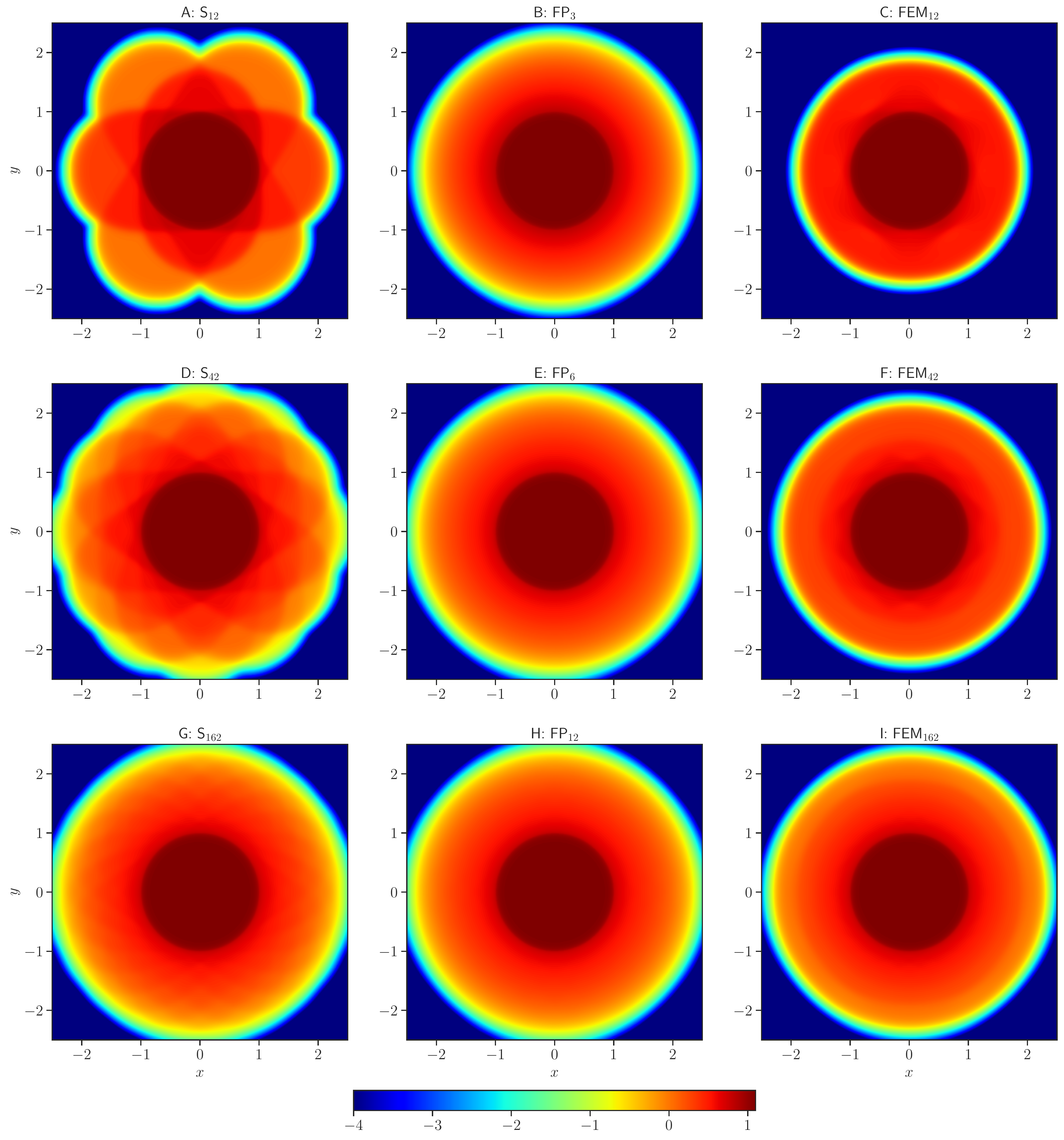}
	\caption{A comparison of~$\log_{10} E$ computed using the~FEM${}_{N}$,~\sn and~\fpn schemes for the homogeneous cylinder problem evolved till $t=1.5$. At this time, radiation reaches the edges of the numerical domain. Prominent ``ray effects" are seen for the~\sn scheme at low to intermediate angular resolutions which become less severe with increasing resolution. The~\fpn and~\femn solutions are of comparable quality at intermediate angular resolutions.}
	\label{fig:cylindert1}
\end{figure*}
Another problem we consider is the lattice test in two dimensions~\cite{Bru2002} which consists of a central square emitting region surrounded by highly scattering and highly absorbing square regions. This problem is designed to test the efficiency of the numerical schemes in complex geometries. The numerical domain is~7 units in length and width and has a central emitting region of unit dimensions having emissivity~$\eta = 1/4 \pi$. The emitting region is represented by the red square in Fig.~\ref{fig:lattice0}. This is surrounded by eleven white squares of unit dimensions which act as absorbing regions with~$\kappa_a = 1$. The surrounding blue region and the central emitting square also has a scattering coefficient~$\kappa_s = 10$. For our numerical simulations, we choose~$\delta x = \delta y = 0.02$ with~$\delta t = 0.0064$. We perform two sets of evolutions with this setup with the~$\textrm{FEM}_N$,~\sn and~\fpn schemes. In the first case, the system is evolved up to~$t = 3.2$. The results are plotted as~$\log_{10}$ of radiation energy density as shown in Fig.~\ref{fig:lattice1}. For the~\fpn simulations, a Lanczos filter with~$\sigma_{\textrm{eff}} \geq 5$ is used to ensure solution positivity.
\begin{figure*}[!t] 
	\centering
	\includegraphics[width=0.8\textwidth]{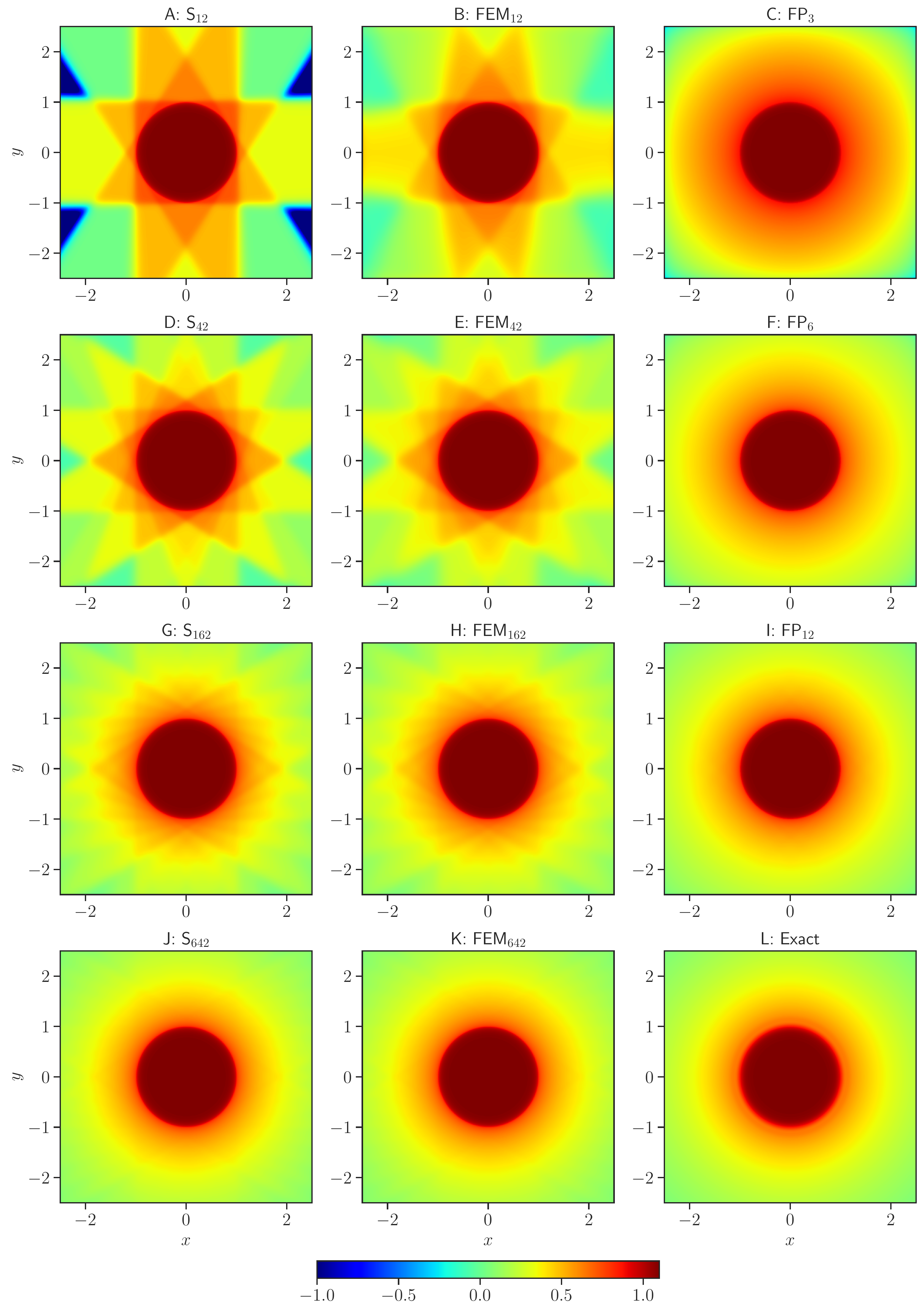}
	\caption{A comparison of steady state values of~$\log_{10} E$ for the homogeneous cylinder test. All runs have been evolved till~$t = 18.75$. Both~\sn and~\femn methods show ray-like artifacts at low angular resolutions with solution quality improving with the number of angles. At~$N = 642$, both~\femn and~\sn solutions are of comparable quality to the exact solution.}
	\label{fig:cylinder2}
\end{figure*}
Around this time, radiation is expected to reach the edges of the boundary, but not at the corners. Except the~\fpn scheme, neither of the other schemes reflect this fact accurately at very low angular resolutions, as can be seen in plots A, B and C of Fig.~\ref{fig:lattice1}. The~$\textrm{S}_{12}$ and~$\textrm{\textit{FEM}}_{12}$ solutions propagate slower than the speed of light with the~\sn method showing significant ray-like artifacts. At higher angular resolutions, solutions for all three schemes tend to attain the correct propagation speed with the~\fpn and~\femn results being of similar quality. Solutions computed with the~\sn method are ``ray-like" even at~$642$ angles and are therefore not preferred for this test. The~\fpn results are dependent on the choice of~$\sigma_{\textrm{eff}}$ with higher opacities resulting in a solution of poorer quality when compared with the~\femn results. Furthermore, these solutions do not propagate at the correct speed. For the second set of runs, the system is evolved up to~$t = 16$ when steady state had been reached. As can be seen in Fig.~\ref{fig:lattice2}, irrespective of the angular resolution chosen, both~\fpn and~\femn schemes give comparable results implying that either method can be chosen for this test. Solutions evaluated by the~\sn method still demonstrate some ray-like artifacts at high angular resolutions and are therefore not preferred over the other two schemes.

\subsection{The homogeneous cylinder test}
\begin{figure*}[!t] 
	\centering
	\includegraphics[width=\textwidth]{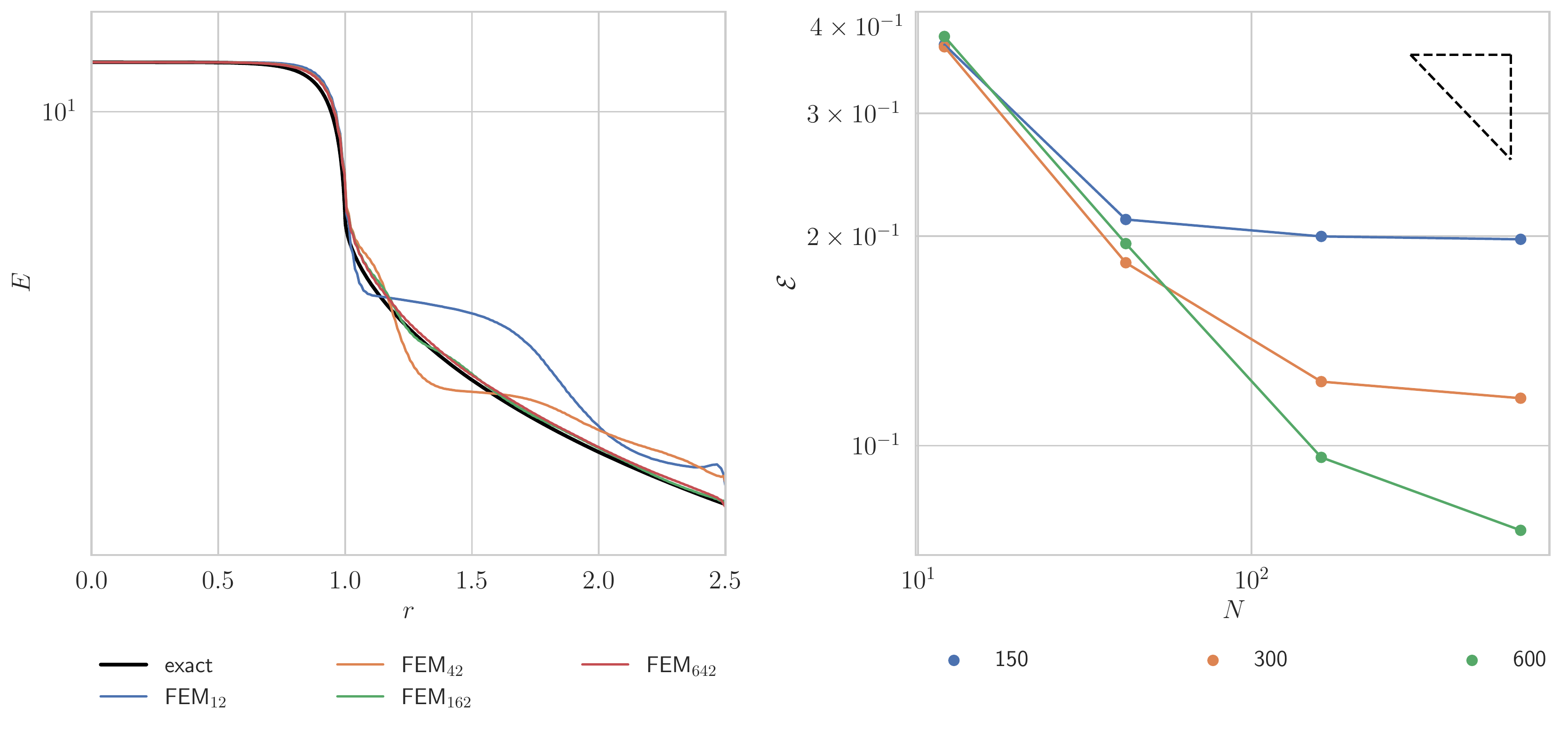}
	\caption{\textit{Left}: A comparison of steady state values of~$E$ along the x-axis between different~\femn resolutions and the exact solution. The~\femn solution is not isotropic at low angular resolutions but approaches the true solution with increasing resolution. \textit{Right}: A comparison of the error~$\mathcal{E}$ between the~\femn at three spatial resolutions. A representative triangle with a hypotenuse of slope $-1/2$ is shown in the upper right corner.}
	\label{fig:cylinder3}
\end{figure*}
The final test we consider is that of an infinite cylindrical source of radiation of unit radius located at the center of the numerical domain. This source has an emissivity~$\eta = 10$ and an absorption coefficient~$\kappa_a = 10$. This problem tests the effectiveness of numerical schemes in handling sharp discontinuities at the surface of the cylindrical source and incorporates the same challenges of the homogeneous sphere problem~\cite{RadAbdRez2013}. The numerical domain extends from~$x,y \in [-2.5,2.5]$ with~$\delta x = \delta y \approx 0.008, 0.017, 0.033$ and~$\delta t = 0.003, 0.0075$. Two sets of simulations are performed, one which is evolved up till~$t = 1.5$ and another up till steady state has been achieved. For the~\fpn runs, ~$\sigma_{\textrm{eff}} = 5$ is chosen as the effective opacity for the Lanczos filer. For the first set of runs, radiation emitted from the cylindrical source should reach the edges of the numerical domain at this time. A comparison of the three schemes is shown in Fig.~\ref{fig:cylindert1}. At low angular resolutions, propagation speeds of solutions in the~\femn and~\sn case are slower than the speed of light. This improves with increase in angular resolution. The~\sn solutions shown ``ray artifacts'' and are of worse quality than the solutions produced by the other two methods. Panels A,B and C of Fig.~\ref{fig:cylindert1} compare solutions at low resolutions. At higher resolutions, as seen in panels D, E, F, G, H and I, the~\femn and~\fpn solutions are of comparable quality. The~\sn solutions agree well with the other methods when~$N \geq 162$.

\begin{figure*}[!t] 
	\centering
	\includegraphics[width=0.5\textwidth]{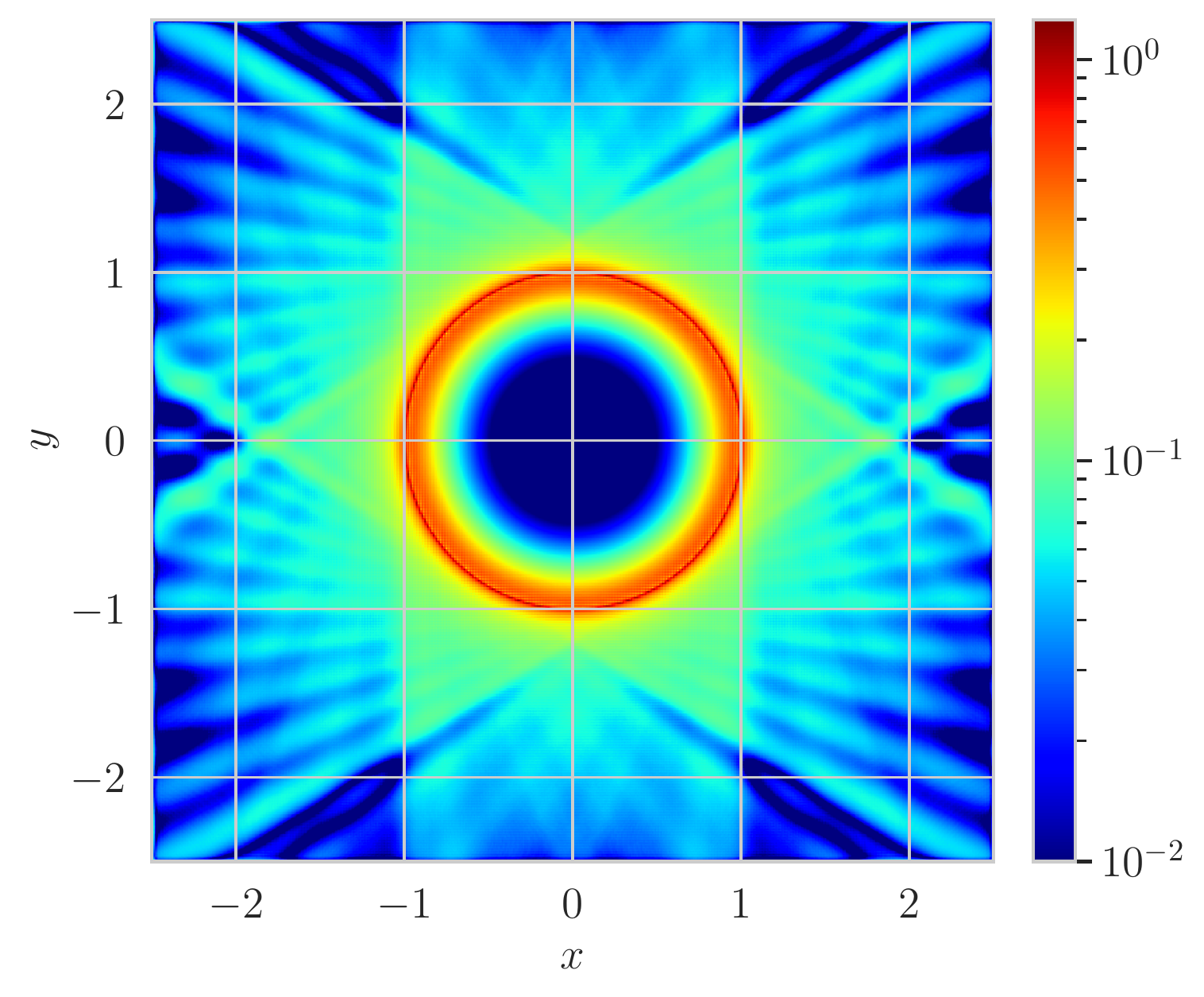}
	\caption{A plot of the absolute value of the difference between the exact solution and the~\femn solution at the highest resolution. The maximum error is seen at the surface of the cylindrical source.}
	\label{fig:cylinder4}
\end{figure*}

For the second set of runs, we evolve the system up till~$t = 18.75$ and compare the results with the exact steady state solution of the problem
\begin{align}
	F (t,r,\phi,\theta) = \frac{B}{\kappa_a} \left(1 - e^{-\kappa_a s}\right), && s(r,\phi,\theta) = \lambda_1(r,\phi,\theta) - \lambda_2(r,\phi,\theta),
\end{align}
where
\begin{align}
	\lambda_1 = \max\left(\frac{r \cos \phi - \sqrt{R - r^2 \sin^2 \phi}}{\sin \theta},0\right), &&
	\lambda_2 = \max\left(\frac{r \cos \phi + \sqrt{R - r^2 \sin^2 \phi}}{\sin \theta},0\right).
\end{align}
Results from this test are shown in Fig.~\ref{fig:cylinder2} where we see ``ray artifacts" for both the~\sn and~\femn solutions at low angular resolutions. We plot the~$L^1$ error in Fig.~\ref{fig:cylinder3} as defined by
\begin{align} \label{eq:cylinder_error}
	\mathcal{E} = \frac{1}{N_p} \sum_{i} \left| E^{\textrm{exact}}_{i} - E^{\textrm{numerical}}_{i} \right|,
\end{align}
where the summation is performed over the entire numerical domain and $N_p$ is the number of points in the spatial domain. At low and intermediate angular resolution, the~\femn solutions are of superior quality to those produced by the other two methods. At very high angular resolutions, the ~\fpn method provide better accuracy. Fig.~\ref{fig:cylinder3} compares the energy density at the highest spatial resolution for four values of $N$. A comparison between the $L^1$ norms at three spatial resolutions~$150$,~$300$ and~$600$ points in each dimension as a function of~$N$ is also shown in the right panel of the same figure. The~\femn method converges to the true solution of the problem with increasing~$N$ at the expected first order provided sufficiently high spatial resolution is chosen for the problem. At low spatial resolutions, however, the error from the spatial discretization dominates as can be seen in Fig.~\ref{fig:cylinder3}. Fig.~\ref{fig:cylinder4} shows that the maximum error is seen at the surface of the cylindrical source.

An advantage of the~\fpn method is the rotational invariance of the scheme. However, in this problem, we find that while~$E$ may remain positive, the distribution function for the~\fpn case may acquire negative values at early times during the simulation. The choice of filter parameters therefore becomes apparent, a problem that is completely avoided by the~\femn or~\sn schemes. We varied the effect opacity of the filter from~$5$ to~$20$ and small negative values still persisted in the solution at early times. Filtering is effective in mitigating this problem only when very high opacity parameters are chosen, which, while ensuring positivity of~$F$ severely degrades the quality of the solution. Due to the higher accuracy at low and intermediate angular resolutions, the~\femn is preferred over the~\sn method unless very high angular resolution is being considered.
\section{Conclusions} \label{sec:conclusions}
We have presented a new numerical scheme for solving the Boltzmann equation using a finite element method in angle where the angular coordinates are discretized using a spherical geodesic grid. This method was then compared with the filtered spherical harmonics scheme~\fpn and the discrete ordinate scheme~\sn using four problems designed to test various aspects of these schemes. Each method has its own advantages and disadvantages. The~\fpn schemes produces good results for problems involving solutions with explicit rotational invariance but require a finely tuned filter parameter to ensure positivity. This choice is problem dependent and is not known beforehand. In these tests, we also find that filtering may still retain small negative values in the solution when the effective filter opacity is kept low while very high values of~$\sigma_{\textrm{eff}}$ generally degrade solution quality at the cost of strict positivity preservation of the distribution function~$F$. Another disadvantage of the~\fpn method is that it performs poorly when problems involve rays or beams of radiation. The~\sn method, on the other hand, ensures that non-physical values do not appear in the solution, but solutions in most cases are contaminated by prominent ``ray effects", especially at low and intermediate angular resolutions.
The~\femn method which we propose in this paper provides an alternative to the~\sn and~\fpn methods with distinct advantages. It overcomes the shortcomings of filtering from the~\fpn scheme by introducing positivity preserving limiters without any arbitrary free parameters. This ensures that the radiation distribution function remains strictly non-negative at all times for all four tests. It also proves to be superior to the~\fpn scheme for handling beams of radiation. Compared to the~\sn scheme, the~\femn scheme proves to be better at mitigating ``ray effects'' and produce significantly superior solutions in all tests except the searchlight test, which yields superior results when the angle of propagation of the beam is a point on the geodesic grid. At high angular resolutions, the~\femn schemes can also handle beams of radiation. The computational cost versus accuracy of the new method, combined with its positivity preservation capabilities and the ability to handle different types of systems makes it a suitable alternative for use in radiation transport problems.

The~\femn method is not without it's disadvantages. For problems with rotational invariance, it produces results which are poorer than the~\fpn method at very low angular resolutions with the appearance of ``ray artifacts''. This can be redressed by increasing angular resolution. Similarly, an accurate treatment of beams with this method demand moderate to high angular resolution. The present method currently deals with the special relativistic scenario without consideration for the energy of radiation carriers. Future work involves addition of frequency dependence in the equations and development of new limiters for ensuring non-negativity in such scenarios. We also intend to extend this treatment to the full general relativistic case.

\section*{Acknowledgments}
\noindent
MKB would like to thank Patrick Mullen for helpful discussions.
We acknowledge funding from the U.S. Department of Energy, Office of
Science, Division of Nuclear Physics under Award Number(s) DE-SC0021177
and from the National Science Foundation under Grants No. PHY-2011725,
PHY-2020275, PHY-2116686, and AST-2108467.
This research used resources of the National Energy Research Scientific
Computing Center, a DOE Office of Science User Facility supported by the
Office of Science of the U.S.~Department of Energy under Contract
No.~DE-AC02-05CH11231. Computations for this research were also
performed on the Pennsylvania State University's Institute for
Computational and Data Sciences' Roar supercomputer.
\bibliographystyle{unsrt}
\bibliography{references}
\end{document}